\documentclass[twocolumn]{autart}    

\usepackage{graphicx}          

\usepackage{amsmath}
\usepackage{amssymb}
\usepackage{amsfonts}
\usepackage{graphicx}
\usepackage{lipsum}
\usepackage{epstopdf}
\usepackage{algorithmic}
\usepackage{hyperref}
\DeclareMathAlphabet{\bit}{OML}{cmm}{b}{it}

\def\<{\leqslant}           
\def\>{\geqslant}           
\def\div{\mathrm{div}}         

\def\d{\partial}

\def\cH{\mathcal{H}}   
\def\mR{\mathbb{R}}    
\def\mC{\mathbb{C}}    

\def\Tr{\mathrm{Tr}}       
\def\rT{\mathrm{T}}        
\def\rF{\mathrm{F}}        

\def\bP{\mathbf{P}}    
\def\bE{\mathbf{E}}    

\def\[[[{[\![\![}
\def\]]]{]\!]\!]}

\def\bra{{\langle}}
\def\ket{{\rangle}}

\def\Bra{\Big\langle}
\def\Ket{\Big\rangle}

\def\re{\mathrm{e}}        
\def\rd{\mathrm{d}}        

\def\cL{\mathcal{L}}

\def\bD{\mathbf{D}}

\def\bF{\mathbf{F}}

\def\x{\times}

\def\fB{\mathfrak{B}}

\def\fF{\mathfrak{F}}

\def\cF{\mathcal{F}}

\def\fW{\mathfrak{W}}

\def\var{\mathbf{var}}

\def\cN{\mathcal{N}}
\def\cS{\mathcal{S}}

\def\mS{\mathbb{S}}

\def\eps{\epsilon}

\begin{document}
\begin{frontmatter}

\title{
Entropy Bounds for
Invariant Measure Perturbations in Stochastic Systems with Uncertain  Noise\thanksref{footnoteinfo}} 

\thanks[footnoteinfo]{The material of this paper was not presented at any conference.}

\author
{Igor G. Vladimirov}\ead{igor.g.vladimirov@gmail.com}

\address
{School of Engineering, Australian National University, Canberra, ACT 2601, Australia}

\begin{abstract}                          
This paper is concerned with stochastic systems whose state is a diffusion process governed by an Ito  stochastic differential equation (SDE).  In the framework of a nominal white-noise model, the SDE is driven by a standard Wiener process. For a scenario of statistical uncertainty, where the driving noise acquires a state-dependent drift and thus deviates from its idealised model, we consider the perturbation of the invariant probability density function (PDF) as a steady-state solution of the Fokker-Planck-Kolmogorov equation. We discuss an upper bound on a logarithmic Dirichlet form for the ratio of the invariant PDF to its nominal counterpart in terms of the Kullback-Leibler relative entropy rate of the actual noise distribution with respect the Wiener measure. This bound  is shown to be achievable,   provided the PDF ratio is preserved by the nominal steady-state probability flux. The logarithmic Dirichlet form bound  is used in order to obtain an upper bound on the relative entropy of the perturbed invariant PDF in terms of quadratic-exponential moments of the noise drift in the uniform ellipticity case. These results are illustrated for perturbations of Gaussian invariant measures in linear stochastic systems involving linear noise drifts.
\end{abstract}


\end{frontmatter}
\endNoHyper

\section{Introduction}

Continuous-time physical systems subject to external random noise 
(and artificial dynamical systems such as stochastic optimization algorithms with intentionally introduced randomness) 
are often described by Ito 
SDEs 
for the evolution of a finite-dimensional state vector driven by a standard Wiener process or more complicated Ito processes \cite{O_2000}. The dynamics of the noise itself  are usually modelled by a shaping filter, also in the form of an SDE,  which can be incorporated in the system, so that the resulting augmented system is driven by a  standard Wiener process. The statistical structure of the system state then corresponds to a diffusion process whose PDF  
(provided it exists and is sufficiently smooth) evolves in time according to the Fokker-Planck-Kolmogorov equation (FPKE)
\cite{BKRS_2015,R_1996,S_2008}  which is a linear second-order 
partial differential equation of parabolic type. 

In the presence of dissipative effects (strong enough to make the system stable), the FPKE usually has a unique normalised steady-state solution describing the invariant measure for the Markovian dynamics of the system state.
However, unmodelled dynamics can lead to deviations of the actual probability law of the driving noise from its  
nominal white-noise model which assumes the noise to be a standard Wiener process. Interpreted as a  statistical uncertainty,  such deviation from the Wiener measure arises, for example, when the external noise acquires a state-dependent drift through a ``parasitic'' feedback coming from the interaction of the system with its environment. The latter  is exemplified by acoustic resonance, 
electromagnetic interference in circuitry, 
or fluid-structure coupling in turbulent flows. 
The presence of an unknown  state-dependent drift in the
noise (making it statistically uncertain and ``coloured'') modifies the drift term of the SDE which governs the system dynamics, and this, in turn,  changes the invariant measure of the system.

Assuming that the FPKE for the perturbed system also has a unique steady-state solution, the present paper investigates the influence of the noise drift 
on the resulting perturbation in the invariant PDF of the system state variables. This perturbation is formulated using the ratio of the invariant PDF for the perturbed system to the  nominal invariant PDF (which the system would have in the nominal white-noise case). We discuss an entropy identity for the noise drift (as a function of the system state) and the logarithmic PDF ratio and derive from it an upper bound for a diffusion-weighted logarithmic Dirichlet form in terms of the second moment of the noise drift over the invariant PDF. Up to a factor of $\frac{1}{2}$,  this moment (which quantifies the ``size'' of the noise drift)
coincides with the rate of the Kullback-Leibler  relative entropy 
\cite{CT_2006} of the probability law of the noise with respect to the Wiener measure.

Entropy and information theoretic  quantification  of statistical uncertainty is employed in a number of approaches to stochastic robust control, including the anisotropy-based theory \cite{VKS_1995} (see also \cite{YKT_2023} and references therein), minimax LQG control  \cite{P_2006,UP_2001} (with its links  \cite{DJP_2000,PUS_2000} to risk-sensitive control),  and other control settings with entropy functionals \cite{CR_2007,TEM_2018}. Entropy theoretic criteria also underlie the Schr\"{o}dinger bridge, which is a state PDF transition problem with relative entropy minimisation, considered both for classical \cite{Beghi_1994,Blaqiere_1992,DaiPra_1991,Mikami_1990} and quantum \cite{BFP_2002} systems  using the formalism of stochastic mechanics \cite{Nelson_2001} (see also  \cite{CL_1994,CGP_2016,PF_2013,VP_2010,VP_2015} and references therein).

Accordingly, the upper bound on the logarithmic Dirichlet form quantifies the influence of the statistical uncertainty in the noise on the deviation of the invariant PDF of the system from its nominal counterpart. We show that this bound is tight in the sense that it is achievable at a particular  noise drift,  provided the PDF ratio is preserved along
the divergenceless probability flux \cite{R_1996} associated with the nominal invariant measure of the system. The logarithmic Dirichlet form  bound  is combined with the logarithmic Sobolev inequality \cite{BE_1985,G_1975,G_1993,L_1992} in order to obtain upper bounds on the Fisher \cite{S_1959} and Kullback-Leibler relative entropies of the invariant PDF of the perturbed system state with respect to its nominal counterpart in the case when the diffusion matrix satisfies the uniform ellipticity condition \cite{E_2008} and the nominal invariant PDF has strong logarithmic  concavity. These results are illustrated for perturbations of Gaussian invariant measures in linear stochastic systems (governed by SDEs with a linear drift and constant diffusion) caused by linear noise drifts.

We mention that regularity of solutions of stationary FPKEs 
was previously studied, for example,  in \cite{BR_1995,BKR_1996,BKR_2006,BKRS_2015,BSV_2016,BKS_2023,MPR_2005}. In particular, their behaviour under drift perturbations for the Ornstein-Uhlenbeck process is considered in \cite{BKS_2023},  and differentiability of the invariant measure over a parameter of the drift vector and diffusion matrix is discussed in \cite{BSV_2016}. Also, Lemma~\ref{lem:iden} of the present paper can be derived by using appropriate modifications of \cite[Theorem 3.1]{BR_1995}, \cite[Theorem 1.1]{BKR_1996},
\cite[Theorem 2.4]{BKR_2006} on logarithmic Dirichlet forms.  However, our results are directly oriented to specific drift perturbations considered in this work under conditions of classical smoothness and vector field decay at infinity. Also, we employ a variational formula from  \cite[Eq. (1.15)]{DE_1997}   (similarly to the way it is used in minimax LQG control)   in order to obtain entropy bounds for the invariant PDF perturbation  in terms of nominal quadratic-exponential moments of the noise drift. Such moments (with integrals or sums of quadratic forms over time) are typical for classical \cite{BV_1985,J_1973,W_1981}  and quantum \cite{B_1996,VPJ_2018,VPJ_2021,V_2022} risk-sensitive control. Furthermore, since the entropy inequalities  use, as  a gain coefficient, the reciprocal to the product of the
uniform ellipticity  and logarithmic concavity constants, we investigate the asymptotic behaviour of the entropy bound for small values of the coefficient and provide a lower bound  for  this coefficient  in terms of the mean divergence of the drift in the nominal system dynamics.

The paper is organised as follows.
Section~\ref{sec:sys} describes the class of stochastic systems with statistically uncertain noise  under consideration.
Section~\ref{sec:inv} specifies the actual and nominal invariant PDFs along with FPKEs for such a system.
Section~\ref{sec:bound} establishes the logarithmic Dirichlet form bound for the PDF ratio.
Section~\ref{sec:flow} discusses the achievability of the upper bound along with the  noise drift which saturates it.
Section~\ref{sec:ent} obtains upper bounds on the Fisher  and Kullback-Leibler relative entropies for the PDF perturbation in the uniform ellipticity case.
Section~\ref{sec:lin} illustrates the results for linear stochastic systems with linear noise drifts  and Gaussian PDFs. Section~\ref{sec:ex}  provides a numerical example on perturbed Langevin dynamics.
Section~\ref{sec:conc} makes concluding remarks.

\section{Uncertain Stochastic Systems}
\label{sec:sys}

We consider a  stochastic system whose state is an $\mR^n$-valued diffusion process $X:=(X_t)_{t\> 0}$
governed by an Ito SDE
\begin{equation}
\label{dX}
    \rd X_t
    =
    f(X_t)\rd t
    +
    g(X_t)
    \rd W_t,
\end{equation}
where the drift vector and the dispersion matrix are specified by given maps $f\in C^1(\mR^n, \mR^n)$ and $g\in C^2(\mR^n, \mR^{n\x m})$ in the appropriate spaces of (respectively, once and twice) continuously differentiable functions. The SDE (\ref{dX}) is
driven by an external random noise $W:=(W_t)_{t\> 0}$ which is an $\mR^m$-valued Ito process on an underlying  probability space $(\Omega, \cF, \bP)$ with a filtration $\fF:= (\cF_t)_{t\> 0}$ (of $\sigma$-algebras $\cF_s \subset \cF_t \subset \cF$ of events for any $t\> s \> 0$)    satisfying the usual conditions \cite{KS_1991}, with the initial system state $X_0$ being $\cF_0$-measurable.   The stochastic differential of $W$ is assumed to be in the form
\begin{equation}
\label{dW}
    \rd W_t
    :=
    h(X_t) \rd t + \rd V_t
\end{equation}
(with zero initial condition $W_0=0$,  without loss of generality),
where the drift vector is described by a map $h\in C^1(\mR^n, \mR^m)$, and $V:= (V_t)_{t\> 0}$ is an $\mR^m$-valued  standard  Wiener process with respect to the filtration $\fF$ (and hence, independent of $X_0$).
In what follows, $\fF$ is assumed to be the natural filtration of the pair $(X_0,V)$, so that for any $t\> 0$, the $\sigma$-algebra $\cF_t$ is generated by $X_0$ and the past history $V_{[0,t]}$ of the process $V$ over the time interval $[0,t]$:
\begin{equation}
\label{cFt}
    \cF_t := \sigma\{X_0, V_{[0,t]}\}.
\end{equation}
The presence of an unknown  drift in (\ref{dW}) makes $W$ a statistically uncertain ``coloured'' noise,  which (in contrast to $V$) is, in general,  correlated with $X_0$. As a consequence, for any fixed but otherwise arbitrary time horizon $T>0$,  the  probability distribution $Q_T = \bP\circ W_{[0,T]}^{-1}$ of the random element $W_{[0,T]}$ (as an  $\cF_T/\fB_T$-measurable map from $\Omega$ to $\fW_T$)    can differ from the Wiener measure $Q_T^*$ on the measurable space $(\fW_T, \fB_T)$. Here, $\fW_T:= C_0([0,T], \mR^m)$ is the  Banach space of continuous functions $w: [0,T]\to \mR^m$ vanishing at the origin and endowed  with the uniform norm,   and $\fB_T$ is the Borel $\sigma$-algebra  generated by open subsets of $\fW_T$. The Wiener measure $Q_T^*$ plays the role of an idealised white-noise model for $W$, which would hold if  the SDE (\ref{dW}) were driftless, leading to $W=V$ and the nominal system dynamics
\begin{equation}
\label{dXnom}
    \rd X_t
    =
    f(X_t) \rd t
    +
    g(X_t)
    \rd V_t.
\end{equation}
The above described  statistical uncertainty (where $Q_T$ deviates from $Q_T^*$) 
results from a feedback loop (due to interaction of the system with its environment) depicted schematically in Fig.~\ref{fig:sys}, which refers to the integral form of the SDEs (\ref{dX}), (\ref{dW}).
\begin{figure}[htbp]
\unitlength=0.65mm
\linethickness{0.2pt}
\centering
\begin{picture}(70,48.50)
    \put(35,45){\makebox(0,0)[cc]{\scriptsize system}}
    \put(-2.5,17.5){\dashbox(45,30)[cc]{}}
    \put(0,35){\framebox(10,10)[cc]{{$\int$}}}
    \put(5,30){\vector(0,1){5}}
    \put(10,40){\vector(1,0){7}}
    \put(20,40){\circle{6}}
    \put(20,40){\makebox(0,0)[cc]{\small$+$}}
    \put(35,40){\vector(-1,0){12}}
    \put(35,40){\line(0,-1){10}}
    \put(0,20){\framebox(10,10)[cc]{{\small $f$}}}
    \put(30,20){\framebox(10,10)[cc]{{\small $g\cdot$}}}
    \put(15,0){\framebox(10,10)[cc]{{$h$}}}
    \put(50,25){\circle{6}}
    \put(50,25){\makebox(0,0)[cc]{\small$+$}}

    \put(45,0){\framebox(10,10)[cc]{{$\int$}}}
    \put(50,10){\vector(0,1){12}}
    \put(20,37){\vector(0,-1){27}}
    \put(20,25){\vector(1,0){10}}
    \put(20,25){\vector(-1,0){10}}
    \put(20,25){\makebox(0,0)[cc]{\tiny$\bullet$}}
    \put(25,5){\vector(1,0){20}}

    \put(47,25){\vector(-1,0){7}}
    \put(65,25){\vector(-1,0){12}}
    \put(69,25){\makebox(0,0)[cc]{\scriptsize $V$}}
    \put(23,28){\makebox(0,0)[cc]{\scriptsize $X$}}
    \put(45,28){\makebox(0,0)[cc]{\scriptsize $W$}}
\end{picture}\vskip-1.5mm
\caption{The stochastic system (\ref{dX}) with an uncertain external noise $W$ arising from a  feedback loop, through which the current system state $X_t$ enters the drift $h(X_t)$ of $W$ in (\ref{dW}), with $(g\cdot W)_t:= \int_0^t g(X_s)\rd W_s = \int_0^t g(X_s) h(X_s)\rd s + (g\cdot V)_t$ using the Ito integral.
}
\label{fig:sys}
\end{figure}

In order to guarantee the existence and uniqueness of a solution  for the SDE
\begin{equation}
\label{dX1}
    \rd X_t
    =
    f_h(X_t)\rd t
    +
    g(X_t)
    \rd V_t,
    \qquad
    f_h:= f+gh,
\end{equation}
obtained by substituting (\ref{dW}) into (\ref{dX}), and also for its nominal counterpart (\ref{dXnom}),
it is assumed (see, for example,  \cite[Theorem~5.2.1]{O_2000}) that the initial state has finite second moments ($\bE(|X_0|^2) <+\infty$, where $\bE(\cdot)$ is expectation),   
and the maps $f$, $g$, $h$, $f_h$ are Lipschitz (and hence,  have at most linear growth at infinity). Under Novikov's condition  \cite{N_1973}
\begin{equation}
\label{Nov}
    \bE_*
    \re^{\frac{1}{2}\int_0^T |h(X_t)|^2\rd t}
    <+\infty,
\end{equation}
where $|\cdot|$ is the usual Euclidean norm, and   $\bE_*(\cdot)$ is the nominal expectation (which considers $W$ in (\ref{dX}) a standard Wiener process independent of $X_0$),  the Kullback-Leibler  relative entropy
of the distribution $Q_T$ of $W_{[0,T]}$ (conditioned on $X_0$)
with respect to
the Wiener measure $Q_T^*$ admits the representation
\begin{equation}
\label{RT}
    R_T
     :=
    \bD(Q_T \| Q_T^*)
    :=
    \bE
        \ln
    \frac{\rd Q_T}{\rd Q_T^*}
    =
        \frac{1}{2}
        \int_0^T
        \bE(|h(X_t)|^2)\rd t.
\end{equation}
Recall \cite{CT_2006}   that the relative entropy 
of a probability measure $M$ on a measurable space $(\fW, \fB)$ with respect to a reference probability measure $N$ on the same space  is defined as $\bD(M\|N):= \bE\ln\frac{\rd M}{\rd N} = \int_\fW \ln\frac{\rd M}{\rd N}(w)M(\rd w)$ (under absolute continuity $M\ll N$, with the logarithm  of the Radon-Nikodym derivative $\frac{\rd M}{\rd N}: \Omega \to \mR_+$ being averaged  over $M$ as the actual probability measure). This definition is combined in (\ref{RT}) with Girsanov's theorem \cite{Girsanov_1960}, 
whereby the corresponding log-likelihood ratio takes the form
\begin{align}
\nonumber
    \ln
    \frac{\rd Q_T}{\rd Q_T^*}
    & =
    \int_0^T
    h(X_t)^\rT\rd W_t  - \frac{1}{2}
    \int_0^T
    |h(X_t)|^2\rd t\\
\label{lnPP}
    & =
    \int_0^T
    h(X_t)^\rT\rd V_t  + \frac{1}{2}
    \int_0^T
    |h(X_t)|^2\rd t,
\end{align}
so that (\ref{RT}) results from averaging (\ref{lnPP}) in accordance with the actual noise and system dynamics (\ref{dW}), (\ref{dX1}), including 
the fact that $V$ is a standard Wiener process independent of $X_0$.
The noise relative entropy $R_T$  in (\ref{RT}) quantifies the deviation of the probability law of  $W$ from that of $V$ (that is, the Wiener measure) over the time interval $[0,T]$ and thus vanishes in the nominal case when $h\equiv 0$.

\section{Actual and Nominal Invariant Measures}
\label{sec:inv}

Under the parabolic H\"{o}rmander condition \cite{H_1967,S_2008}, the value $X_t$  of the diffusion process $X$ in (\ref{dX1})  at any time $t>0$ is an absolutely continuously distributed  random vector with a smooth PDF $p_t: \mR^n \to \mR_+$ (with respect to the Lebesgue measure in $\mR^n$) satisfying the FPKE
\begin{equation}
\label{FPKE}
    \d_t p_t
    =
    \cL^{\dagger}(p_t),
\end{equation}
where $\cL^\dagger$ is a linear differential operator acting on a function $\varphi\in C^2(\mR^n, \mR)$ as
\begin{equation}
\label{cL+}
    \cL^{\dagger}(\varphi)
    :=
    \frac{1}{2}
    \div^2(\varphi D)
    -\div (\varphi f_h),
\end{equation}
and the map $D \in C^2(\mR^n, \mS_n^+)$,  with values in the set $\mS_n^+ \subset \mS_n$ of positive semi-definite matrices in the subspace $\mS_n$ of real symmetric matrices of order $n$,  describes the diffusion matrix:
\begin{equation}
\label{D}
    D(x) := g(x)g(x)^{\rT},
    \qquad
    x:=(x_k)_{1\< k\<n } \in \mR^n.
\end{equation}
Here, repeated application of the divergence operator to a function $G :=
(G_{jk})_{1\<j,k\< n} \in C^2(\mR^n,\mS_n)$ yields
the maps
$    \div G
    :=
    (
        \sum_{k=1}^{n}
        \nabla_k G_{jk}
    )_{1\< j\< n} \in C^1(\mR^n,
\mR^n)
$
and $
    \div^2 G
    :=
    \div \div G
    =
    \sum_{j,k=1}^n
    \nabla_j \nabla_k G_{jk} \in C(\mR^n, \mR)
$,
where $\nabla_k:= \d_{x_k}$ is the partial derivative over the $k$th
Cartesian coordinate in $\mR^n$. In accordance with the 
divergence theorem, 
$\cL^{\dagger}$ in (\ref{cL+}) is the formal adjoint \cite{S_2008} of the
infinitesimal generator $\cL$ for the diffusion process $X$, acting on a function $\varphi \in C^2(\mR^n, \mR)$ with the gradient vector $\nabla \varphi: = (\nabla_k \varphi)_{1\< k\< n}$ and the Hessian matrix $\varphi'': = (\nabla_j \nabla_k \varphi)_{1\< j,k\< n}$ as
\begin{equation}
\label{cL}
    \cL(\varphi)
    :=
    f_h^{\rT}\nabla \varphi
    +
    \frac{1}{2}
    \bra D,  \varphi'' \ket_\rF,
\end{equation}
where $\bra M, N \ket_\rF := \Tr (M^\rT N)$ is the Frobenius inner product \cite{HJ_2007} of real matrices. Using the notation $\bra \alpha, \beta \ket:= \int_{\mR^n} \alpha(x)\beta(x)\rd x$ for functions $\alpha, \beta: \mR^n\to \mR$ with an integrable point-wise product $\alpha\beta \in L^1(\mR^n, \mR)$, a sufficient condition for the adjointness relation
$\bra u, \cL(v) \ket = \bra \cL^\dagger(u), v \ket$ to hold is provided by the property that  the pair $(u, v)$ is $(f_h,D)$-\emph{decaying}  as specified below.

\begin{defn}
\label{def:auv}
For given maps $\phi\in C^1(\mR^n,\mR^n)$ and $G \in C^2(\mR^n, \mS_n)$, the pair $(u,v)$ of functions  $u,v \in C^2(\mR^n, \mR)$ is said to be $(\phi,G)$-decaying if the following vector field (as an $\mR^n$-valued function of $x \in \mR^n$) satisfies
\begin{align}
\nonumber
    a_{u,v}
    & :=
    \Big(
        \phi - \frac{1}{2}\div G
    \Big)
    uv +\frac{1}{2}
    G(u\nabla v - v\nabla u)\\
\label{auv}
    & =
    o(|x|^{1-n}),
    \qquad
    {\rm  as}\
    x \to \infty.
\end{align}
\end{defn}
The decay condition in (\ref{auv}) (which is weaker than the requirement that $u$ or $v$ has compact support)   is a joint property of $\phi$, $G$, $u$, $v$ at infinity. Its meaning is clarified by integrating both sides of the identity $
    (\phi^\rT \nabla v + \frac{1}{2}\bra G, v''\ket_\rF)u
     -
    (\frac{1}{2}\div^2(uG)- \div(u\phi)
    )v = \div a_{u,v}
$ (similar to those in Green's formulas) over a ball $\{x \in \mR^n: |x|\< \rho\}$ of radius $\rho$    and applying the divergence theorem along
with the fact that the decay rate of the vector field $a_{u,v}$  in (\ref{auv}) makes its flux through the sphere $\cS_\rho := \{x \in \mR^n: |x|=\rho\}$ asymptotically  vanish:  $\lim_{\rho\to +\infty} \int_{\cS_\rho} a_{u,v}\cdot \rd \sigma=0$, where $\rd \sigma$ is the oriented  surface area element (with the outward normal). For any function $s \in C^2(\mR^n, \mR)$ with nonzero values, $a_{u/s, s v} = a_{u,v} + uv G\nabla \ln |s|$.     In the case of $uv\ne 0$, this transformation follows from the representation
\begin{equation}
\label{auv0}
    a_{u,v}   =
    \Big(
        \phi
    +\frac{1}{2}
    \Big(
        G\nabla \ln \Big|\frac{v}{u}\Big|
        -\div G
    \Big)
    \Big) uv.
\end{equation}

We assume that the FPKE (\ref{FPKE}) has a unique steady-state solution $p \in C^2(\mR^n, \mR_+)$ (with the normalization $\bra 1, p\ket=1$):
\begin{equation}
\label{p}
  \cL^\dagger(p) = 0.
\end{equation}
The invariant PDF $p$ is, in general,  different from  another invariant PDF $p_*$ which the system would have in the nominal case of $h\equiv 0$ in (\ref{dW}). The nominal invariant PDF $p_*$, which is also assumed to exist and be unique in $C^2(\mR^n, \mR_+)$ subject to the normalization $\bra 1, p_*\ket =1$,  satisfies
\begin{equation}
\label{p*}
  \cL_*^\dagger(p_*) = 0,
\end{equation}
where
\begin{align}
\label{cL*+}
    \cL_*^{\dagger}(\varphi)
    & :=
    \frac{1}{2}
    \div^2(\varphi D)
    -
    \div (\varphi f),\\
\label{cL*}
    \cL_*(\varphi)
    & :=
    f^{\rT}\nabla \varphi
    +
    \frac{1}{2}
    \bra D,  \varphi''\ket_\rF
\end{align}
are the nominal versions of (\ref{cL+}), (\ref{cL}), respectively,  with $h\equiv 0$.

\section{Logarithmic Dirichlet Form Bound}
\label{sec:bound}

In  what follows, we will discuss an entropy identity in order to obtain an upper bound for the discrepancy between the invariant PDF $p$ of the perturbed system and its nominal counterpart $p_*$. The bound will be formulated in terms of a steady-state second moment of the noise drift map $h$ over $p$. To this end, we will use the PDF ratio
\begin{equation}
\label{r}
  r(x):= \frac{p(x)}{p_*(x)},
  \qquad
  x \in \mR^n,
\end{equation}
which is well-defined on the set $\{x \in \mR^n: p(x)>0\}$,  provided $p$ vanishes whenever $p_*$ does, thus making the invariant measure of the system absolutely continuous with respect to the nominal invariant  measure. For simplicity, both PDFs $p$ and $p_*$ are assumed to be positive everywhere, and hence, the property  $r \in C^2(\mR^n, \mR_+)$ is inherited from them.   Also, we will use a function $\psi\in C^1(\mR^n,\mR^m)$ associated by
\begin{equation}
\label{psi}
  \psi(x):= g(x)^\rT \nabla \ln r(x)
\end{equation}
with the dispersion matrix map $g$ and the PDF ratio (\ref{r}). In view of the factorisation (\ref{D}) of the diffusion matrix,
\begin{equation}
\label{ggD}
    |\psi| = \|\nabla \ln r\|_D,
\end{equation}
where  $\|v\|_M:= |\sqrt{M}v| = \sqrt{v^{\rT}M v}$ is a weighted Euclidean (semi-) norm of a real vector $v$, specified by a  real positive (semi-) definite symmetric  matrix $M$. If $p=p_*$ everywhere in $\mR^n$, then the ratio (\ref{r}) reduces to $r\equiv 1$, and the function $\psi$ in (\ref{psi}) (and its point-wise norm (\ref{ggD})) is the identical zero.

\begin{lem}
\label{lem:iden}
Suppose the  generator $\cL$ in (\ref{cL}),  evaluated at the logarithmic PDF ratio $\ln r$ from (\ref{r}), the   noise drift map $h$ in (\ref{dW})  and the  function $\psi$ in (\ref{psi}) have the following finite moments
\begin{equation}
\label{moms}
    \bE|\cL(\ln r)|<+\infty,\
  \bE (|h|^2)<+\infty,\
  \bE (|\psi|^2) < +\infty
\end{equation}
over the  invariant PDF $p$ of the perturbed system (\ref{dX1}) in (\ref{p}).
Also, suppose the pair $(p,\ln r)$ is $(f_h,D)$-de\-cay\-ing, while $(p_*, r)$ is $(f,D)$-decaying in the sense of Definition~\ref{def:auv}.
Then
\begin{equation}
\label{id0}
    \bE\Big(h^\rT \psi  - \frac{1}{2} |\psi|^2\Big)
    = 0.
\end{equation}
\end{lem}
\begin{pf}
In application to the generator $\cL$ in (\ref{cL}) and the PDF ratio $r$ in  (\ref{r}),
the Fleming logarithmic transformation \cite{F_1982}
(see also \cite[Eq.~(81) on p.~201]{BFP_2002})  leads to
\begin{align}
\nonumber
    \cL(\ln r)
     & =
     \frac{1}{r}
    \cL(r)
    -
    \frac{1}{2}
    \|\nabla \ln r\|_D^2\\
\label{logtrans}
    & =
     \frac{1}{r}
    (\cL_*(r) + \Lambda(r))
    -
    \frac{1}{2}
    |\psi|^2.
\end{align}
Here, use is also made of (\ref{ggD}) and an auxiliary  differential operator $\Lambda$,  which acts on a function $\varphi\in C^1(\mR^n, \mR)$ as
\begin{equation}
\label{Lam}
    \Lambda(\varphi)
    :=
    (\cL-\cL_*)(\varphi)
    =
    (f_h-f)^\rT \nabla \varphi
    =
    h^\rT g^\rT \nabla \varphi,
\end{equation}
in accordance with (\ref{dX1}), (\ref{cL}), (\ref{cL*}). Application of (\ref{Lam}) to $\varphi :=r$  yields
\begin{equation}
\label{Lam/r}
    \frac{1}{r} \Lambda(r)
    =
    \frac{1}{r} h^\rT g^\rT \nabla r
    =
    h^\rT g^\rT \nabla \ln r
    =
    h^\rT \psi,
\end{equation}
where the last equality uses (\ref{psi}). By combining (\ref{logtrans}) with (\ref{Lam/r}), it follows that
\begin{equation}
\label{LL}
    \frac{1}{r}\cL_*(r)
    =
    \cL(\ln r) - h^\rT \psi + \frac{1}{2}|\psi|^2.
\end{equation}
In view of  (\ref{moms}) and the Cauchy-Bunyakovsky-Schwarz inequality
\begin{equation}
\label{ineq1}
    |\bE (h^\rT \psi)| \< \sqrt{\bE (|h|^2) \bE (|\psi|^2)},
\end{equation}
the right-hand side of (\ref{LL}) has a finite expectation $\bE(\cdot)$ over the PDF $p$, and hence, so also does its left-hand side:
\begin{equation}
\label{ELL}
    \bE
    \Big(\frac{1}{r}\cL_*(r)\Big)
    =
    \bE\cL(\ln r) - \bE \Big(h^\rT \psi - \frac{1}{2}|\psi|^2\Big).
\end{equation}
Now, since the pair $(p,\ln r)$ is assumed to be $(f_h,D)$-de\-cay\-ing, then
\begin{equation}
\label{ELlnr}
    \bE\cL(\ln r)
    =
    \bra
        p,
        \cL(\ln r)
    \ket
    =
    \bra
        \cL^\dagger(p),
        \ln r
    \ket
    =
    0,
\end{equation}
where the rightmost equality follows from (\ref{p}). Due to the PDF change identity $\frac{p}{r} = p_*$ from (\ref{r}), the left-hand side of (\ref{ELL}) can  be represented in terms of the nominal expectation $\bE_*(\cdot)$:
\begin{align}
\nonumber
    \bE
    \Big(\frac{1}{r}\cL_*(r)\Big)
    & =
    \Bra
        p,
         \frac{1}{r}
        \cL_*(r)
    \Ket
    =
    \overbrace{\bra p_*, \cL_*(r)\ket}^{\bE_*\cL_*(r)}\\
\label{Lpr}
    & =
    \bra \cL_*^\dagger(p_*), r \ket        = 0,
\end{align}
where the last two equalities use (\ref{p*}) and the condition that the pair $(p_*, r)$ is $(f,D)$-decaying. The relation (\ref{id0}) is now obtained by substitution of (\ref{ELlnr}),  (\ref{Lpr}) into (\ref{ELL}). 
\hfill$\blacksquare$
\end{pf}

Note that the identity (\ref{id0}) is of entropy theoretic nature. More precisely,
suppose the system (\ref{dX1}) is initialised at its invariant PDF $p$ from (\ref{p}). Then  its state vector $X_t$ retains this PDF at any time $t\>0$. In particular, assuming that $\bE |\ln r|<+\infty$,
the following quantity remains constant:
\begin{equation}
\label{entpp}
    \bE \ln r(X_t)
    =
    \int_{\mR^n} p(x) \ln r(x) \rd x
    =
    \bD(p\| p_*),
\end{equation}
which, in view of (\ref{r}),  is the relative entropy of the invariant PDF $p$ with respect to its nominal counterpart $p_*$ (we slightly abuse notation by referring to the densities instead of the distributions). On the other hand, from
the Ito lemma \cite{KS_1991} combined with (\ref{dX1}), (\ref{cL}),  (\ref{psi})  along with the first and third conditions in (\ref{moms}), it follows that $\ln r(X_t)$ is an Ito process with the stochastic differential
\begin{align}
\label{dlnr}
    \rd \ln r(X_t)
    = &
    \cL(\ln r)(X_t)\rd t + \rd Z_t.
\end{align}
Here,
\begin{equation}
\label{Mt}
    Z_t:= \int_0^t \psi(X_t)^\rT \rd V_t
\end{equation}
is a square integrable martingale \cite{O_2000} (with zero mean $\bE Z_t = 0$ and variance $\bE (Z_t^2) =  t\bE (|\psi|^2)$ for any $t\> 0$ by the Ito isometry)   with respect to the natural filtration $\fF$ of the pair $(X_0, V)$ given by (\ref{cFt}).   Therefore, in view of (\ref{dlnr}),
 the equalities (\ref{ELlnr}) represent the fact that the time derivative of the relative entropy
(\ref{entpp}) vanishes:
$
    (\bE \ln r(X_t))^{^\centerdot}
     =
    \bE \cL (\ln r)(X_t) = 0
$, which is not affected by the martingale (\ref{Mt}).

The following theorem is a direct corollary of Lemma~\ref{lem:iden}.
\begin{thm}
\label{th:dev}
Under the conditions of Lemma~\ref{lem:iden},
the mean square of the function $\psi$ in 
(\ref{psi}) 
is bounded in terms of that of the noise drift map $h$ in (\ref{dW}) as
\begin{equation}
\label{dev}
    \bE (|\psi|^2)
       \<
      4 \bE (|h|^2).
\end{equation}
\end{thm}
\begin{pf}
By combining (\ref{id0}) with (\ref{ineq1}), it follows that
\begin{equation}
\label{ineq2}
    (\bE (|\psi|^2))^2
    =
    (2
    \bE (h^\rT \psi))^2
    \<
    4
    \bE (|h|^2) \bE (|\psi|^2),
\end{equation}
and hence,
$
    \bE (|\psi|^2)
    \<
    4
    \bE (|h|^2)
$, which establishes (\ref{dev}). An alternative (yet equivalent) argument employs the determinant of a positive semi-definite Gram matrix 
as
\begin{align}
\nonumber
    0
    & \<
    \det
    \bE
    \begin{bmatrix}
        |\psi|^2 & h^\rT \psi\\
        h^\rT \psi & |h|^2
    \end{bmatrix} =
    \bE (|\psi|^2)\bE (|h|^2) - (\bE (h^\rT \psi))^2\\
\label{ineq3}
    & =
    \bE (|\psi|^2)
    \Big(
        \bE (|h|^2) - \frac{1}{4}\bE (|\psi|^2)
    \Big),
\end{align}
which also uses the relation $\bE (h^\rT \psi) = \frac{1}{2}\bE(|\psi|^2)$ from  (\ref{id0}) and leads to (\ref{dev}).
\hfill$\blacksquare$
\end{pf}

Recalling (\ref{ggD}), note that the nonnegative quantity $\bE(|\psi|^2) =  \bE (\|\nabla \ln r\|_D^2)$  on the left-hand side of (\ref{dev}) is a diffusion-weighted Dirichlet form \cite{E_2008} for the logarithmic PDF ratio.

In view of (\ref{RT}) (assuming that the sufficient condition (\ref{Nov}) is satisfied), the inequality (\ref{dev}) can also be represented as
\begin{equation}
\label{Epsi2}
    \bE(|\psi|^2)
    \<
    8 \lim_{T\to +\infty}
    \frac{R_T}{T} ,
\end{equation}
where the limit is the relative entropy rate of the input noise and coincides with $\frac{R_T}{T} = \frac{1}{2}\bE(|h|^2)$ for any time horizon $T>0$ if the initial system  state $X_0$ has the invariant PDF $p$. This limit quantifies the ease of detecting the drift in the noise $W$  based on the strong law of large numbers which manifests itself in sufficiently long  samples of the log-likelihood ratio (\ref{lnPP}) (see, for example, \cite{LS_2001}). The latter consideration motivates the following definition.

\begin{defn}
The noise drift map $h$ in (\ref{dW}) is  said to be $\gamma$-stealthy  if its mean square over the invariant PDF $p$ satisfies
\begin{equation}
\label{gam}
    \bE (|h|^2) \< 2 \gamma,
\end{equation}
where $\gamma\> 0$ is a given threshold on the relative entropy rate in (\ref{RT}).
\end{defn}

\section{Upper Bound Achievability}
\label{sec:flow}

The inequality  (\ref{dev}) is tight in the sense that it becomes an equality (as do the inequalities in (\ref{ineq1}), (\ref{ineq2}), (\ref{ineq3})) if the noise drift map $h$ in (\ref{dW}) is related to the PDF ratio $r$  through  (\ref{psi}) by
\begin{equation}
\label{hpsi}
    h = \frac{1}{2}\psi.
\end{equation}
However, this relation is not necessarily compatible with $p$ being an invariant PDF satisfying (\ref{p}) and $p_*$ its nominal counterpart in (\ref{p*}). A condition for such compatibility  is provided by the following theorem in terms of an auxiliary vector field (a steady-state probability flux \cite{R_1996}) in $\mR^n$ given by
\begin{equation}
\label{vel}
        U:=
        p_*f
        -
        \frac{1}{2}
        \div(p_* D),
\end{equation}
which is divergenceless since $\div U = -\cL_*^\dagger(p_*) = 0$ in view of (\ref{p*}), (\ref{cL*+}).

\begin{thm}
\label{th:comp}
Suppose the noise drift map $h$ is related by (\ref{hpsi}) to the PDF ratio $r$ in  (\ref{r}) through $\psi$ in  (\ref{psi}).  Then
\begin{equation}
\label{Ur}
  U(x)^\rT \nabla r(x) = 0,
  \qquad
  x \in \mR^n,
\end{equation}
that is,
$r$ is preserved by the incompressible flow generated by the vector field $U$ in (\ref{vel}).
\end{thm}
\begin{pf}
By using the identity $p = r p_*$ from (\ref{r}), it follows that
\begin{align}
\label{divpf}
  \div(pf)
  = &
  \div(r p_* f)
  =
  r \div(p_* f) + p_* f^\rT \nabla r,\\
\nonumber
    \div^2(pD)
    = &
    \div^2(r p_* D)
    =
    \div(r\div(p_* D) + p_* D\nabla r)\\
\label{div2pD}
     =&
    r\div^2(p_* D)
    +
    \div(p_* D)^\rT \nabla r
    +
    \div(p_* D\nabla r).
\end{align}
Therefore, substitution of  (\ref{divpf}), (\ref{div2pD}) into (\ref{cL*+}) leads to
\begin{align}
\nonumber
    \cL_*^\dagger(p)
    = &
    \frac{1}{2}
    \div^2(p D)
    -\div (p f)\\
\nonumber
    = &
    \frac{1}{2}
    \Big(
            r\div^2(p_* D)
    +
    \div(p_* D)^\rT \nabla r
    +\div(p_* D\nabla r)
    \Big)\\
\nonumber
    & -r \div(p_* f) - p_* f^\rT \nabla r\\
\nonumber
    = &
        r \cL_*^\dagger(p_*)
        -
        U^\rT
        \nabla r+
        \frac{1}{2}
        \div(p_* D \nabla r)\\
\label{cL*+p}
    = &
        \frac{1}{2}
        \div(p_* D \nabla r)
        -
        U^\rT
        \nabla r,
\end{align}
where use is made of (\ref{p*}) along with the vector field $U$ from  (\ref{vel}). Now,
if the relation (\ref{hpsi}) is satisfied, then $gh = \frac{1}{2} gg^\rT \nabla \ln r = \frac{1}{2r} D \nabla r$ in view of (\ref{D}), and hence, the adjoint of the operator $\Lambda$ in (\ref{Lam}) acts on the invariant PDF $p$  as
\begin{equation}
\label{Lam+p}
    \Lambda^\dagger(p)
     =
    (\cL^\dagger-\cL_*^\dagger)(p)
    =
    -\div(pgh)
    =
    -
    \frac{1}{2}
    \div(p_* D \nabla r),
\end{equation}
By taking the sum of the right-hand sides of (\ref{cL*+p}) and  (\ref{Lam+p}), it follows that in the case of the noise drift map (\ref{hpsi}),
\begin{equation}
\label{cL+p}
  \cL^\dagger(p)
  =
  \cL_*^\dagger(p) + \Lambda^\dagger(p)
  =
    -U^\rT
  \nabla r.
\end{equation}
A comparison of (\ref{cL+p}) with (\ref{p}) leads to (\ref{Ur}), which is equivalent to $U$ in (\ref{vel}) being tangent to the level  hypersurfaces of the PDF ratio $r$, or the property  that  $r$ is a first integral of motion for the incompressible flow in $\mR^n$ generated by $U$. 
\hfill$\blacksquare$
\end{pf}

The preservation of the PDF ratio $r$ along the vector field $U$ in Theorem~\ref{th:comp} is thus necessary for achievability of the upper bound (\ref{dev}), and the noise drift (\ref{hpsi}) saturates it.

\section{Entropy Bounds on PDF perturbation}
\label{sec:ent}

The inequality (\ref{dev}) is particularly useful 
when the diffusion matrix map $D$ in (\ref{D}) satisfies 
the uniform ellipticity condition \cite{E_2008}
\begin{equation}
\label{lam}
    \lambda
    :=
    \inf_{x \in \mR^n}
    \lambda_{\min}(D(x))>0
\end{equation}
(with $\lambda_{\min}(\cdot)$ the smallest eigenvalue of a real symmetric matrix),  and hence, with necessity, the dispersion matrix $g(x)$ is of full row rank for all $x \in \mR^n$, so that $n\< m$.   In this case, (\ref{ggD})  admits a lower bound $|\psi| \> \sqrt{\lambda} |\nabla \ln r|$ which  leads to the inequality
\begin{equation}
\label{bF}
    \bF(p \| p_*)
     :=
    \bE_*
    \Big(
        \frac{1}{r}
        |\nabla r|^2
    \Big)
    =
    \bE (|\nabla \ln r|^2)
    \<
    \frac{1}{\lambda}\bE (|\psi|^2),
\end{equation}
whose left-hand side is the Fisher relative entropy  (originating from the Fisher information for the shift parameter; see, for example,  \cite[Eq. (1.4)]{S_1959}). By the probabilistic version \cite{BE_1985} of the  logarithmic Sobolev inequalities \cite{G_1975,G_1993,L_1992}, the Kullback-Leibler relative entropy $\bD(p\| p_*)$ in (\ref{entpp}) satisfies
\begin{equation}
\label{bDbF}
    \bD(p\| p_*)
    =
    \bE_* (r \ln r)
    \<
    \frac{1}{2\mu} \bF(p \| p_*),
\end{equation}
provided the function $\ln p_*$ 
is strongly concave in the sense that its Hessian matrix is uniformly negative definite:
\begin{equation}
\label{mu}
    \mu
    :=
    -
    \sup_{x\in \mR^n}
    \lambda_{\max}
    ((\ln p_*)''(x))
    >0,
\end{equation}
where $\lambda_{\max}(\cdot)$ is the largest eigenvalue.  In contrast to $\lambda$ in (\ref{lam}), the quantity  $\mu$ in (\ref{mu}) depends not only on the diffusion part of the nominal SDE  (\ref{dXnom}) for the system, but also on  the drift map $f$ (we will address this dependence  in Lemma~\ref{lem:K0}).
A combination of (\ref{bF}) with (\ref{bDbF}) and (\ref{dev}) of Theorem~\ref{th:dev} leads to
\begin{equation}
\label{bD}
    \bD(p\| p_*)
     \<
    \frac{K}{2} \bE(|\psi|^2)
    \<
    2K \bE(|h|^2),
\end{equation}
where
\begin{equation}
\label{K}
    K := \frac{1}{\lambda \mu}
\end{equation}
is an auxiliary quantity which has the physical dimension of time.
In particular,  for $\gamma$-stealthy noise drift maps $h$ (as specified by the relative entropy  threshold in (\ref{gam}) and motivated by (\ref{Epsi2}) under the sufficient condition (\ref{Nov})), the inequality (\ref{bD}) yields
\begin{equation}
\label{bDgam}
    \bD(p\| p_*)
    \<
    4K \gamma.
\end{equation}
The last inequality describes the influence of such deviations of the actual noise $W$ in (\ref{dW}) from its nominal white-noise  model on the discrepancy between the corresponding PDFs $p$ and  $p_*$, so that, up to a factor of $4$, the quantity    $K$ from (\ref{K}) plays the role of a gain coefficient  for this influence.

While the condition (\ref{gam}) and the corresponding bound (\ref{bDgam}) are formulated in terms of $\bE(|h|^2)$ over the invariant PDF $p$,
it is also possible to develop (\ref{bD}) in an alternative direction, which (instead of (\ref{gam})) employs moments of $h$ over the nominal invariant PDF $p_*$ and takes into account the deviation of $p$ from $p_*$ in computing the right-hand side of (\ref{bD}).
 More precisely, it follows from (\ref{bD}) that
\begin{equation}\
\label{bD1}
    \bD(p\| p_*)
     \<
    2K \Phi(\bD(p\| p_*)),
\end{equation}
where $\Phi: \mR_+ \to \mR_+$ is a nondecreasing function defined as
\begin{equation}
\label{Phi}
    \Phi(\eps)
    :=
    \sup_{p:\, \bD(p\|p_*)\< \eps}
    \bE(|h|^2),
    \qquad
    \eps\> 0,
\end{equation}
that is,  the largest mean square of $h$ over those PDFs $p$ whose deviation from $p_*$ is bounded in terms of the relative entropy constraint.  In particular, $\Phi(\eps) \> \Phi(0)= \bE_* (|h|^2)$. The maximal solution of (\ref{bD1}) provides the upper bound
\begin{equation}
\label{bDupper}
    \bD(p\|p_*)
    \<
    \sup
    \{
        \eps:\ 0\<  \eps \< 2K \Phi(\eps)
    \}.
\end{equation}
Now, the function $\Phi$ in (\ref{Phi}) can be computed through the variational formula of \cite[Proposition 1.4.2 on pp. 33, 34]{DE_1997} as
\begin{equation}
\label{dual}
    \Phi(\eps)
    =
    \inf_{\theta \>0}
        \frac{\Psi(\theta) + \eps}{\theta},
\end{equation}
where the ratio is set to its limit values $\Psi'(0)$ or $+\infty$ by continuity   at $\theta=0$   if $\eps=0$ or $\eps>0$, respectively. Here, use is made of
the cumulant-generating  function (CGF) $\Psi$ for $|h|^2$ associated by
\begin{equation}
\label{Psi}
  \Psi(\theta)
  :=
  \ln \Xi(\theta),
  \qquad
  \theta \>0,
\end{equation}
with a quadratic-exponential moment of the noise drift map $h$ over the PDF $p_*$:
\begin{equation}
\label{Xi}
    \Xi(\theta)
    :=
    \bE_* \re^{\theta |h|^2}
    =
    \bra
        p_*,
        \re^{\theta |h|^2}
    \ket,
\end{equation}
that is, the nominal moment-generating function for $|h|^2$. Here, $\theta$ plays the role of a risk-sensitivity parameter and, similarly to $K$, has the physical dimension of time.
With $\Xi(0)=1$, we assume that $\Xi(\theta)$ is finite at least for some $\theta>0$, so that
\begin{equation}
\label{theta*}
    \theta_*
    :=
    \sup
    \{
        \theta\> 0:\
        \Xi(\theta)< +\infty
    \}>0,
\end{equation}
and hence, $\Xi(\theta)\> 1$ (accordingly, $\Psi(\theta)\> 0$)
for any $\theta \in [0, \theta_*)$. Note that if the system state is initialised at the PDF $p_*$ and obeys the nominal SDE (\ref{dXnom}), then (\ref{theta*}) secures  the fulfillment of the Novikov condition (\ref{Nov}) for any $0<T <2\theta_*$:
\begin{align*}
    \bE_* & \re^{\frac{1}{2}\int_0^T |h(X_t)|^2 \rd t}
     =
    \bE_* \re^{\frac{T}{2}\frac{1}{T}\int_0^T |h(X_t)|^2 \rd t}\\
    & \<
    \frac{1}{T}
    \int_0^T
    \bE_* \re^{\frac{T}{2}|h(X_t)|^2}
    \rd t =
    \Xi(T/2) < +\infty,
\end{align*}
where Jensen's inequality is combined with the convexity of the exponential function, and (\ref{Xi}) is used.
Except for the trivial case of $|h|$ being an identical  constant (we exclude it from consideration),   the CGF $\Psi$ in (\ref{Psi})  is strictly increasing and strictly convex, with positive first and second derivatives $\Psi'(\theta)$, $\Psi''(\theta)$,   which,  together with
$
    \d_\theta\frac{\Psi(\theta) + \eps}{\theta}
    =
    \frac{\nu(\theta)-\eps}{\theta^2},
$
implies that
for any $\eps\> 0$, the infimum in (\ref{dual}) is achieved at a unique point $\theta = \nu^{-1}(\eps)$. Here, $\nu^{-1}$ is the inverse of  a strictly increasing differentiable  function   $\nu: [0, \theta_*) \to \mR_+$ defined as the Bregman divergence \cite{B_1967}  for $\Psi$ at the points $0$, $\theta$ (with $\Psi(0)=0$):
\begin{equation}
\label{nu}
    \nu(\theta)
    :=
    \theta\Psi'(\theta)-\Psi(\theta),
    \qquad
    0\< \theta < \theta_*,
\end{equation}
so that, due to the representation (\ref{dual}), the function $\Phi$ in (\ref{Phi}) takes the form
\begin{equation}
\label{PhiPsi'}
    \Phi(\eps)
    =
    \frac{\Psi(\nu^{-1}(\eps)) + \eps}{\nu^{-1}(\eps)}
    =
    \Psi'(\nu^{-1}(\eps)).
\end{equation}
Since the strictly increasing function $\nu^{-1}$, associated with (\ref{nu}) and used in (\ref{PhiPsi'}),  maps $\mR_+$ onto the interval $[0, \theta_*)$, then the set on the right-hand side of (\ref{bDupper}) is the image
\begin{equation}
\label{chiTheta}
    \{
        \eps:\ 0\<  \eps \< 2K \Phi(\eps)
    \}
    =
    \nu(\Theta_K)
\end{equation}
of another set
\begin{equation}
\label{Theta}
    \Theta_K
    :=
    \{
    \theta \in [0, \theta_*):\
    \nu(\theta) \< 2K \Psi'(\theta)
    \}.
\end{equation}
By the structure (\ref{nu}) of the function $\nu$, the inequality in (\ref{Theta}) is equivalent to
\begin{align}
\nonumber
    \nu_K(\theta)
    & :=
    \nu(\theta) - 2K \Psi'(\theta)\\
\label{nuK}
    & =
    (\theta-2K)\Psi'(\theta)-\Psi(\theta) \< 0,
\end{align}
whose left-hand side satisfies $\nu_K(2K) = -\Psi(2K)<0$ and is strictly increasing over $\theta>2K$ since
\begin{equation}
\label{nuder}
    \d_\theta \nu_K(\theta) = (\theta-2K)\Psi''(\theta) >0,
    \qquad
    \theta \in (2K, \theta_*),
\end{equation}
where
$\Psi''(\theta) >0$ in accordance  with the strict convexity of $\Psi$ mentioned above. 
Hence, if
\begin{equation}
\label{Ksmall}
  K < \frac{1}{2}\theta_*,
\end{equation}
then the maximal solution $\theta_K$ of the inequality (\ref{nuK}), which is related to the set (\ref{Theta}) by $\theta_K= \sup \Theta_K$, is found uniquely from
\begin{equation}
\label{sol0}
    \nu_K(\theta_K) = 0
\end{equation}
and satisfies
\begin{equation}
\label{theta2K}
 2K < \theta_K < \theta_*.
\end{equation}
Therefore, due to (\ref{chiTheta}), the upper bound (\ref{bDupper})    acquires  the following  form (its proof is provided by the above discussion).

\begin{thm}
\label{th:bDnu}
Suppose the uniform ellipticity (\ref{lam}) for the diffusion matrix and the strong logarithmic concavity (\ref{mu}) for the nominal invariant PDF $p_*$ are satisfied. Also, suppose the function $|h|: \mR^n\to \mR_+$ for the noise drift map $h$ in (\ref{dW}) is not identically constant and satisfies (\ref{theta*}).  Then the relative entropy in (\ref{entpp}) admits the upper bound
\begin{equation}
\label{bDupper1}
    \bD(p\|p_*)
    \<
    \nu(\theta_K) = 2K \Psi'(\theta_K)
\end{equation}
in terms of the gain coefficient $K$ from  (\ref{K}) and the unique solution $\theta_K$ of (\ref{sol0}), where the function $\nu_K$ is associated by (\ref{nuK}) with the nominal CGF $\Psi$ of $|h|^2$ from (\ref{Psi}).
\end{thm}

The inequality (\ref{bDupper1}) is also applicable to other measures of deviation of $p$ from $p_*$,  which are bounded in terms of the Kullback-Leibler relative entropy. For example, a combination of Pinsker's inequality \cite[Lemma 2.5 on p. 88]{T_2009} with (\ref{bDupper1}) leads to an upper bound for the $L^1$-distance:
\begin{equation}
\label{L1}
    \int_{\mR^n} |p(x)-p_*(x)|\rd x
    \<
    \sqrt{2\bD(p\|p_*)}
    \<
    \sqrt{2 \nu(\theta_K)}.
\end{equation}
The relations (\ref{Ksmall})--(\ref{bDupper1}) (and their corollary (\ref{L1}))  can be reformulated in a ``small-gain'' fashion: for any fixed but otherwise arbitrary $\eps>0$, the following  implication holds:
\begin{equation}
\label{Keps}
    K \< \frac{\eps}{2\Psi'(\nu^{-1}(\eps))}\
    \Longrightarrow\
    \bD(p\| p_*) \< \eps.
\end{equation}
In turn, the entropy bound (\ref{Keps}) can be applied not only to other distance measures for the PDFs as in (\ref{L1}), but also to obtaining guaranteed upper bounds on a cost functional in the form of a generalised moment $\bE \phi = \bra p, \phi\ket$ of the system variables specified by a function $\phi: \mR^n\to \mR_+$:
$
    \sup_{p:\, \bD(p\|p_*)\< \eps}
    \bE\phi
    =
    \inf_{\vartheta  \>0}
        \frac{\ln \bE_*\re^{\vartheta\phi} + \eps}{\vartheta}
$,
which,  similarly to (\ref{Phi}), (\ref{dual})--(\ref{Xi}), employs the variational formula together with the nominal CGF $\vartheta\mapsto \ln \bE_*\re^{\vartheta\phi}$ for $\phi$.

In practical calculations, the solution of the equation (\ref{sol0}) for computing $\nu(\theta_K)$ as a function of $K$, or   the inversion of $\nu$ for finding the critical value of  $K$ versus $\eps$ in (\ref{Keps}),  can be avoided by parameterising this curve in the $(K,\eps)$-plane as
\begin{equation}
\label{tpar}
    [0, \theta_*)
    \ni \theta \mapsto
    \Big(\frac{\nu(\theta)}{2\Psi'(\theta)}, \nu(\theta)\Big).
\end{equation}
For small values of $K$, 
the upper bound in (\ref{bDupper1}) behaves asymptotically as $\nu(\theta_K) \sim 2K \bE_*(|h|^2)$, which corresponds to neglecting the discrepancy between $p$ and $p_*$ on the right-hand side of
(\ref{bD}).
The next term in the asymptotic expansion of $\nu(\theta_K)$, as $K\to 0+$, which takes into account the deviation of $p$ from $p_*$, is provided below.

\begin{thm}
\label{th:asy}
Suppose $|h|$ 
is not identically constant,  and (\ref{theta*}) is satisfied.
Then the first two leading terms in the asymptotic expansion of the upper bound (\ref{bDupper1}) are as follows:
\begin{align}
\nonumber
    \nu(\theta_K)
    = &
    2 \bE_*(|h|^2) K
    +
    4 \sqrt{\bE_*(|h|^2) \var_*(|h|^2)}\, K^{3/2}\\
\label{nuasyK}
    &  + o(K^{3/2}),
    \qquad
    {\rm as}\ K \to 0+,
\end{align}
where
\begin{equation}
\label{var*}
    \var_*(|h|^2) =
    \bE_*(|h|^4) - (\bE_*(|h|^2))^2
\end{equation}
is the variance of $|h|^2$ over the nominal invariant PDF $p_*$.
\end{thm}
\begin{pf}
By the implicit function theorem, the solution $\theta_K$ of  (\ref{sol0})  is a strictly increasing differentiable function of $K\in (0, \frac{1}{2}\theta_*)$,   with
$\rd \theta_K/\rd K = -\frac{\d_K \nu_K(\theta)}{\d_\theta \nu_K(\theta)}\big|_{\theta = \theta_K}>0$ in view of (\ref{nuder}), (\ref{theta2K}) and since $\d_K \nu_K(\theta) = -2\Psi'(\theta)<0$ from (\ref{nuK}) for any $\theta \in [0,\theta_*)$. This  monotonicity implies the existence of a limit $\theta_0:=  \lim_{K\to 0+}\theta_K\> 0$. From the convergence $\nu_K(\theta)\to \nu(\theta)$, as $K\to 0+$, which is uniform over $\theta\in [0, \delta]$ for any $\delta\in [0, \theta_*)$, it follows that $\nu(\theta_0) = \lim_{K\to 0+}\nu_K(\theta_K) = 0$ due to (\ref{sol0}), and hence,  $\theta_0 = 0$:
\begin{equation}
\label{theta0}
    \lim_{K\to 0+}\theta_K = 0.
\end{equation}
By (\ref{Psi})--(\ref{nu}), the functions $\nu$, $\Psi'$  satisfy the asymptotic relations
\begin{align}
\label{nuasy}
    \nu(\theta)
    & =
    \theta^2
    \Big(
        \frac{\Psi(\theta)}{\theta}
    \Big)'
    =
    \frac{1}{2}
    \Psi''(0)
    \theta^2
    +
    O(\theta^3) ,\\
\label{Psi'asy}
    \Psi'(\theta)
    & =
    \Psi'(0) + \theta \Psi''(0)
    +
    O(\theta^2),
    \quad
    {\rm as}\
    \theta \to 0.
\end{align}
Here, the property $\Psi(0)=0$ is used along with the first two cumulants of $|h|^2$ over the nominal invariant PDF $p_*$:
\begin{equation}
\label{Psi'Psi''}
    \Psi'(0)
     =
    \bE_*(|h|^2),
    \qquad
    \Psi''(0)
     = \var_*(|h|^2),
\end{equation}
where the second cumulant is the nominal variance of $|h|^2$ given by (\ref{var*}). Due to (\ref{theta0}), a combination of (\ref{nuasy}), (\ref{Psi'asy}) with the equality on the right-hand side of   (\ref{bDupper1}) leads to the asymptotic equivalence
$
    \frac{1}{2}\Psi''(0)\theta_K^2 \sim 2\Psi'(0) K$,
and hence,
\begin{equation}
\label{tasy}
    \theta_K \sim 2\sqrt{\frac{\Psi'(0)}{\Psi''(0)}K},
    \qquad
    {\rm as}\
    K \to 0+.
\end{equation}
Therefore, $O(\theta_K^2) = O(K) = o(\sqrt{K})$, and substitution of (\ref{tasy}) into (\ref{Psi'asy}) yields
\begin{equation}
\label{Psi'tK}
    \Psi'(\theta_K)
    =
    \Psi'(0) + 2\sqrt{\Psi'(0)\Psi''(0)K} + o(\sqrt{K}).
\end{equation}
The relation (\ref{nuasyK}) for the right-hand side of (\ref{bDupper1}) now follows from (\ref{Psi'tK}) in view of (\ref{Psi'Psi''}). \hfill$\blacksquare$
\end{pf}

Since the entropy bounds (\ref{bDgam}), (\ref{bDupper1}) involve the coefficient $K$ from (\ref{K}), and the asymptotic relation (\ref{nuasyK}) pertains to its small values, of relevance is the following lower bound for $K$ in terms of the drift part of the nominal system dynamics (\ref{dXnom}).

\begin{lem}
\label{lem:K0}
Suppose the uniform ellipticity (\ref{lam}) and strong logarithmic  concavity (\ref{mu}) conditions are satisfied. Also, suppose
\begin{equation}
\label{moms*}
  \bE_*|\cL_*(\ln p_*)| < +\infty,
  \quad
  \bE_*|\div f| < +\infty,
\end{equation}
and $(p_*, \ln p_*)$ is $(f,D)$-decaying, while $(1,p_*)$ is $(f,0)$-decaying  (see Definition~\ref{def:auv}).
Then $K$ in (\ref{K})  satisfies
\begin{equation}
\label{Kdiv}
  K \> -\frac{n}{2\bE_* \div f}.
\end{equation}
\end{lem}
\begin{pf}
Similarly to (\ref{ELlnr}), (\ref{Lpr}) in  the proof of Lemma~\ref{lem:iden},  from the first condition in (\ref{moms*}), the  assumption that $(p_*, \ln p_*)$ is $(f,D)$-decaying and from (\ref{p*}), it follows that
\begin{equation}
\label{EEE}
    \bE_* \cL_*(\ln p_*)
     =
    \bra p_*, \cL_*(\ln p_*)\ket
     =
    \bra \cL_*^\dagger(p_*), \ln p_*\ket
    = 0.
\end{equation}
The second condition in (\ref{moms*}) and the assumption that $(1,p_*)$ is $(f,0)$-decaying (that is, $p_*(x)f(x) = o(|x|^{1-n})$, as $x\to \infty$, in accordance with (\ref{auv}) with $\phi:=f$, $G:=0$, $u:=1$, $v:=p_*$) imply that
\begin{align}
\nonumber
    \bE_*\div f
    & =
    \bra p_*, \div f\ket
    =
    -\bra 1, f^\rT \nabla p_* \ket\\
\label{E*div}
    & =
    -\bra p_*, f^\rT \nabla \ln p_* \ket
    =
    -
    \bE_* (f^\rT \nabla \ln p_*).
\end{align}
At the same time, (\ref{cL*}) yields $\cL_*(\ln p_*)=f^\rT \nabla \ln p_* + \frac{1}{2}\bra D, (\ln p_*)''\ket_\rF$,  whose averaging over $p_*$, in combination with  (\ref{EEE}), (\ref{E*div}),  leads to
\begin{equation}
\label{E*div1}
    \bE_*\div f
    =
    \frac{1}{2}
    \bE_* \bra D, (\ln p_*)''\ket_\rF.
\end{equation}
By the monotonicity of the Frobenius inner product $\bra \cdot, \cdot\ket_\rF$ on the cone $\mS_n^+$, the inequalities (\ref{lam}), (\ref{mu}) imply that
\begin{align}
\nonumber
    \bra D(x), &(\ln p_*)''(x)\ket_\rF
    =
    -\bra D(x), -(\ln p_*)''(x)\ket_\rF    \\
\nonumber
    & \<
    -\lambda \bra I_n, -(\ln p_*)''(x)\ket_\rF
    =
    \lambda \Tr (\ln p_*)''(x)\\
\label{Dlnp}
    & \< -\lambda \mu n = -\frac{n}{K},
    \qquad
    x \in \mR^n,
\end{align}
where the last equality uses (\ref{K}).  A combination of (\ref{E*div1}) with (\ref{Dlnp}) leads to $\bE_* \div f\< -\frac{n}{2K}<0$, which establishes (\ref{Kdiv}).
\hfill$\blacksquare$
\end{pf}

The averaging in (\ref{Kdiv}) is redundant if the map $f$ is affine, when $\div  f$ is identically  constant,  as is the case for linear stochastic systems.

\section{Illustration for Linear-Gaussian Dynamics}
\label{sec:lin}

Consider a  class of linear stochastic systems described by (\ref{dX}) with a linear drift vector and a constant dispersion matrix:
\begin{equation}
\label{AB}
  f(x) := Ax,
  \quad
  g(x) := B,
  \qquad
  x \in \mR^n,
\end{equation}
where $A \in \mR^{n\x n}$, $B \in \mR^{n\x m}$ are given matrices, so that the corresponding diffusion matrix in (\ref{D}) is also constant:
\begin{equation}
\label{DBB}
  D = BB^\rT.
\end{equation}
Assuming that $A$ is Hurwitz, the unperturbed system (\ref{dXnom}), governed by the SDE
\begin{equation}
\label{Xlin*}
    \rd X_t = AX_t \rd t + B\rd V_t
\end{equation}
(with the standard Wiener process $V$),
has a unique nominal invariant Gaussian measure $\cN(0, P_*)$ with zero mean and covariance matrix
\begin{equation}
\label{P*}
    P_*
    =
    \int_{\mR_+} \re^{tA} D\re^{tA^\rT} \rd t
    =
    \frac{1}{2\pi}
    \int_{\mR}
    \Sigma(\omega)
    \rd\omega.
\end{equation}
Here,
\begin{equation}
\label{Sigma}
      \Sigma(\omega)
      :=
      F(i\omega) F(i\omega)^*,
      \qquad
      \omega \in\mR
\end{equation}
(with $(\cdot)^* := (\overline{(\cdot)})^{\rT}$  the complex conjugate transpose) is the  nominal spectral density of the stationary Gaussian process $X$, and
\begin{equation}
\label{F}
    F(s):= (sI_n - A)^{-1}B,
    \qquad
    s \in \mC,
\end{equation}
 is the transfer function from the incremented  input $V$ to $X$.
The matrix $P_*$ in (\ref{P*}) is the infinite-horizon controllability Gramian \cite{KS_1972} of the pair $(A,B)$ satisfying the algebraic Lyapunov equation (ALE)
\begin{equation}
\label{P*ALE}
    AP_* + P_*A^\rT + D = 0
\end{equation}
in view of (\ref{DBB}).
If $(A,B)$ is controllable, then $P_*\succ 0$, and the nominal invariant measure is absolutely continuous with the PDF
\begin{equation}
\label{p*gauss}
    p_*(x)
    =
    \frac{(2\pi)^{-n/2}}{\sqrt{\det P_*}}
    \re^{-
        \frac{1}{2}
        \|x\|_{P_*^{-1}}^2},
    \qquad
    x \in \mR^n.
\end{equation}
In the presence of a linear drift in the actual noise dynamics (\ref{dW}), given by
\begin{equation}
\label{hlin}
    h(x) := Nx,
    \qquad
    x \in \mR^n,
\end{equation}
with a constant matrix $N \in \mR^{m\x n}$ (specifying the parasitic coupling of the noise to the system through the feedback loop in Fig.~\ref{fig:sys}), the nominal SDE (\ref{Xlin*}) is replaced with
\begin{equation}
\label{Xlin}
    \rd X_t = (A + BN)X_t \rd t + B\rd V_t,
\end{equation}
in accordance with (\ref{dX1}). Assuming that $A+BN$ is also Hurwitz,  the system state has a zero-mean Gaussian invariant measure $\cN(0,P)$ whose covariance matrix $P$ satisfies the ALE
\begin{equation}
\label{PALE}
    (A+BN)P + P(A+BN)^\rT + D = 0.
\end{equation}
Since the pair $(A+BN,B)$ inherits controllability from $(A,B)$, then $P$ is also positive definite, giving rise to the invariant PDF
\begin{equation}
\label{pgauss}
    p(x)
    =
    \frac{(2\pi)^{-n/2}}{\sqrt{\det P}}
    \re^{-
        \frac{1}{2}
        \|x\|_{P^{-1}}^2},
    \qquad
    x \in \mR^n.
\end{equation}
The deviation of the invariant covariance matrix $P$ from its nominal counterpart   $P_*$  depends on the matrix $N$  and vanishes when $N=0$. The following theorem, which specifies Lemma~\ref{lem:iden} and Theorem~\ref{th:dev} in the case of linear-Gaussian dynamics, provides an identity and an upper bound for the difference
\begin{equation}
\label{PP}
    \Pi:= P_*^{-1}-P^{-1} = \Pi^\rT
\end{equation}
of the corresponding precision matrices in terms of the Frobenius norm $\|M\|_\rF := \sqrt{\bra M, M\ket_\rF}$ for real matrices. Their applicability is secured by the fact that the corresponding vector fields in (\ref{auv}) (see also (\ref{auv0})) are organised as a product of a Gaussian PDF and a polynomial (with an exponentially fast decay at infinity) and since the integrability conditions (\ref{moms}) hold for a similar reason.

\begin{thm}
\label{th:devlin}
Suppose the pair $(A,B)$ in (\ref{AB}) is controllable, and $A$ is Hurwitz. Then for any matrix $N$ in (\ref{hlin}) such that $A+BN$ is also Hurwitz, the precision matrix difference (\ref{PP}) for the nominal (\ref{Xlin*}) and perturbed (\ref{Xlin})  system dynamics,  satisfies
\begin{align}
\label{PPprod}
    \Bra
        BN - \frac{1}{2} D\Pi,
        \Pi P
    \Ket_\rF & = 0,\\
\label{PPnorm}
    \|\sqrt{D} \Pi \sqrt{P}\|_\rF
      & \<  2 \|N\sqrt{P}\|_\rF.
\end{align}
\end{thm}
\begin{pf}
In the linear-Gaussian setting under consideration, substitution of (\ref{AB}), (\ref{p*gauss}), (\ref{pgauss}) into (\ref{r}) leads to
\begin{equation}
\label{lnr'}
  \nabla \ln r(x)  = \Pi x,
  \qquad
  x \in \mR^n,
\end{equation}
with the matrix $\Pi$ from (\ref{PP}), so that the map $\psi$  in (\ref{psi}) is also linear:
\begin{equation}
\label{psilin}
    \psi(x) = B^\rT \Pi x.
\end{equation}
Hence, the identity (\ref{PPprod}) is obtained by computing the expectations of quadratic forms in (\ref{id0}) over the invariant Gaussian distribution $\cN(0,P)$ as
\begin{align}
\nonumber
    0
    & = \bE\Big(h^\rT \psi  - \frac{1}{2} |\psi|^2\Big)
    =
    \Tr
    \Big(
        \Big(
            N^\rT B^\rT \Pi - \frac{1}{2} \Pi D \Pi
        \Big)P
    \Big)\\
\label{id0lin}
    & =
    \Bra
        BN - \frac{1}{2} D\Pi,
        \Pi P
    \Ket_\rF,
\end{align}
where (\ref{DBB}), (\ref{hlin}), (\ref{psilin}) are used along with the symmetry of the matrix $\Pi$ in (\ref{PP}). In  a similar fashion,
the averaging in (\ref{dev}) yields
\begin{align}
\nonumber
    \bE(|\psi|^2)
    & =
    \bra
     \Pi D\Pi, P
    \ket_\rF
    =
    \|\sqrt{D} \Pi \sqrt{P}\|_\rF^2
    \<
    4 \bE(|h|^2)\\
\label{devlin}
    &
    =
    4
    \bra
        N^\rT N, P
    \ket_\rF
    = 4 \|N\sqrt{P}\|_\rF^2,
\end{align}
which establishes the norm bound (\ref{PPnorm}) (whose derivation from (\ref{PPprod}) by matrix manipulations without the probabilistic reasoning behind (\ref{id0}), (\ref{dev}),  (\ref{id0lin}), (\ref{devlin}) would be less intuitive).
\hfill$\blacksquare$
\end{pf}

We will now apply Theorem~\ref{th:comp} to achievability of an equality in (\ref{PPnorm}).
Since, in the linear-Gaussian case,   the diffusion matrix $D$ in (\ref{DBB}) is constant, then $\div(p_* D) = D\nabla p_*$, and hence, the vector field $U$ in  (\ref{vel}) takes the form
\begin{align}
\nonumber
        U(x)
        & =
        p_*(x)
        \Big(
        f(x)
        -
        \frac{1}{2}
        D \nabla \ln p_*(x)
        \Big)\\
\label{Ulin}
        & =
        p_*(x)
        H x,
        \qquad
        x \in \mR^n.
\end{align}
Here, use is also made of (\ref{AB}), (\ref{p*gauss}) together with an auxiliary matrix
\begin{equation}
\label{H}
    H:=
        A +
        \frac{1}{2}
        D P_*^{-1}
        =
        \mho  P_*^{-1},
\end{equation}
which is defined in terms of a real antisymmetric matrix
\begin{equation}
\label{mho}
        \mho
        :=
        AP_* +
        \frac{1}{2}
        D = -\mho ^\rT,
\end{equation}
whose antisymmetry follows from the nominal ALE (\ref{P*ALE}) as $\mho  + \mho ^\rT = AP_* + P_* A^\rT + D = 0$.  Note that the matrix $H$ in (\ref{H}) is Hamiltonian in the sense of the symplectic structure specified by $\mho $ (provided the matrix $\mho $ in (\ref{mho}) is nonsingular, in which case, the state dimension $n$ is, with necessity, even). It now follows from (\ref{lnr'}), (\ref{Ulin}), (\ref{r}) that
\begin{align*}
        U(x)^\rT \nabla r(x)
        & =
        p_*(x)
        r(x)
        (\nabla \ln r(x))^\rT
        H x\\
        & =
        p(x)
        x^\rT \Pi
        H x,
        \qquad
        x \in \mR^n,
\end{align*}
with $\Pi$ from (\ref{PP}),
and hence, (\ref{Ur}) is equivalent to the fulfillment of
$
    x^\rT \Pi H x = 0
$ for all $x \in \mR^n$, which holds if and only if the matrix
\begin{equation}
\label{PiH0}
    \Pi H
     =
    P_*^{-1} H -P^{-1}HPP^{-1}
\end{equation}
of this quadratic form is antisymmetric. In turn,  the latter condition reduces to the antisymmetry  of the matrix $HP$ (since $P_*^{-1}H = P_*^{-1}\mho  P_*^{-1}$ in (\ref{PiH0}) is antisymmetric):
\begin{equation}
\label{PiH}
    HP + PH^\rT = 0.
\end{equation}
Therefore, by comparing (\ref{hpsi}) with (\ref{psilin}) and using  Theorem~\ref{th:comp}, it follows that  the inequality (\ref{PPnorm}) becomes an equality if $N = \frac{1}{2} B^\rT \Pi$ in (\ref{hlin}), provided the condition (\ref{PiH}) is satisfied.

We will now apply the entropy bounds of Section~\ref{sec:ent} to the linear-Gaussian setting.  The uniform ellipticity (\ref{lam}) is equivalent to positive definiteness of the diffusion matrix $D$  in (\ref{DBB}):
\begin{equation}
\label{lamgauss}
  \lambda = \lambda_{\min}(D) >0
\end{equation}
(that is, $n\< m$, and $B$ in (\ref{AB}) is of full row rank). This implies the controllability of the pair $(A,B)$, which (together with $A$ being Hurwitz)  makes the nominal invariant covariance matrix $P_*$ positive definite. Since the Hessian matrix $-(\ln p_*)'' = P_*^{-1}$ for the Gaussian PDF $p_*$ in  (\ref{p*gauss}) is also constant, the condition (\ref{mu}) is satisfied:
\begin{equation}
\label{mugauss}
    \mu = \lambda_{\min}(P_*^{-1}) >0.
\end{equation}
The quadratic-exponential moment (\ref{Xi}) of the linear map $h$ from (\ref{hlin}) is calculated (in a standard fashion) as
\begin{align}
\nonumber
    \Xi(\theta)
    & =
    \int_{\mR^n}
    p_*(x)
    \re^{\theta |Nx|^2}
    \rd x\\
\nonumber
    & =
    (2\pi)^{-n/2}
    \int_{\mR^n}
    \re^{\theta |N\sqrt{P_*} u|^2 - \frac{1}{2}|u|^2}
    \rd u\\
\nonumber
    & =
    (2\pi)^{-n/2}
    \int_{\mR^n}
    \re^{-\frac{1}{2}\|u\|_{I_n - 2\theta \sqrt{P_*} N^\rT N\sqrt{P_*}}^2}
    \rd u\\
\nonumber
    & =
    1\Big/\sqrt{\det(I_m - 2\theta N P_* N^\rT)}
\end{align}
(where $I_s$ denotes the identity matrix of order $s$), and
hence, the nominal CGF for $|h|^2$ in (\ref{Psi}) acquires the form
\begin{align}
\nonumber
  \Psi(\theta)
  & =
  -\frac{1}{2}
  \ln\det
  (I_m - 2\theta N P_* N^\rT)\\
\label{Psigauss}
  & =
  \frac{1}{2}
  \sum_{k\> 1}
  \frac{1}{k}
  (2\theta)^k
  \|N\sqrt{P_*}\|_{2k}^{2k}
\end{align}
for any
\begin{equation}
\label{thetamax}
  \theta
  <
  \theta_*
  =
  \frac{1}{2 \|N\sqrt{P_*}\|^2},
\end{equation}
in accordance with (\ref{theta*}).
Here, $\|M\|:= \sqrt{\lambda_{\max}(M^\rT M)} = \lim_{k\to +\infty}\|M\|_{2k}$ is the operator norm (the largest singular value), and $\|M\|_{2k}:= \sqrt[2k]{\Tr((M^\rT M)^k)}$ is the Schatten $2k$-norm \cite[p. 441]{HJ_2007} of a real matrix $M$ (with the Frobenius norm being its particular case: $\|\cdot\|_\rF = \|\cdot\|_2$).  Under the condition (\ref{thetamax}),  the matrix $I_m - 2\theta N P_* N^\rT$ is positive definite, and the derivative of (\ref{Psigauss}) is given by
\begin{align}
\nonumber
  \Psi'(\theta)
  & =
  \bra
    (I_m - 2\theta N P_* N^\rT)^{-1},
    N P_* N^\rT
  \ket_\rF\\
\label{Psi'gauss}
    & =
    \|(I_m - 2\theta N P_* N^\rT)^{-1/2}
    N \sqrt{P_*}\|_\rF^2
\end{align}
(in view of the identity $(\ln \det G)' = \Tr (G^{-1} G')$ for a differentiable nonsingular  matrix-valued  function $G$ of a scalar parameter)
and is strictly positive whenever $N\ne 0$. By (\ref{Psigauss}), (\ref{Psi'gauss}), the corresponding function (\ref{nu}) takes the form
\begin{align}
\nonumber
    \nu(\theta)
    = &
    \frac{1}{2}
    \Big(
        \Tr
        ((I_m - 2\theta N P_* N^\rT)^{-1})
        -m\\
\label{nugauss}
    & +
  \ln\det
  (I_m - 2\theta N P_* N^\rT)
  \Big),
\end{align}
which can be used for computing the upper bound (\ref{bDupper1}) in the elliptic linear-Gaussian case. The coefficients of the Taylor series in (\ref{Psigauss}) provide the cumulants (\ref{Psi'Psi''}) for the asymptotic relation (\ref{nuasyK}):
\begin{equation}
\label{Psi'Psi''gauss}
    \Psi'(0) = \|N\sqrt{P_*}\|_\rF^2,
    \qquad
    \Psi''(0) = 2\|N\sqrt{P_*}\|_4^4.
\end{equation}
The following lemma (which specifies Lemma~\ref{lem:K0} for the linear case) provides an inequality for  the influence of the drift part of the nominal SDE (\ref{Xlin*}) on the coefficient $K$ from (\ref{K}), which is used in the entropy bounds (\ref{bDgam}), (\ref{bDupper1}).

\begin{lem}
\label{lem:K}
Suppose the matrix $A$ in (\ref{AB}) is Hurwitz, and the diffusion matrix $D$ in (\ref{DBB}) is positive definite. Then the coefficient $K$ in (\ref{K}), computed in terms of (\ref{lamgauss}), (\ref{mugauss}), satisfies
\begin{equation}
\label{Klower}
  K \> -\frac{n}{2\Tr A}.
\end{equation}
\end{lem}
\begin{pf}
Although (\ref{Klower}) is a corollary of (\ref{Kdiv}) due to the linearity of the map $f$ in (\ref{AB}), whereby $\div f = \Tr A$, we will also provide an alternative proof using the Hamiltonian structure of the matrix $H$ in (\ref{H}). From the latter property, or directly from  the orthogonality of the subspaces of real symmetric and real antisymmetric matrices in the sense of the Frobenius inner product $\bra\cdot, \cdot \ket_\rF$, it follows   that $\Tr H=0$ (which is closely related to the Liouville theorem \cite{A_1989} on the phase-space volume preservation by Hamiltonian flows). Hence, (\ref{H}) implies that
\begin{equation}
\label{P*ALE2}
    -2 \Tr A
     = \bra D,  P_*^{-1}\ket_\rF
    \>
      \lambda \Tr (P_*^{-1})
      \>
    n \lambda\mu = \frac{n}{K},
\end{equation}
where, similarly to (\ref{Dlnp}), the monotonicity of $\bra \cdot, \cdot\ket_\rF$ on $\mS_n^+$ is combined with $D \succcurlyeq \lambda I_n$, $\Tr (P_*^{-1})\> n \mu$ from (\ref{lamgauss}), (\ref{mugauss}),   and  (\ref{K}) is used. The inequality (\ref{P*ALE2}) leads to  (\ref{Klower}) (note that $\Tr A <0$ follows directly from $A$ being Hurwitz).
\hfill$\blacksquare$
\end{pf}

From the lower bound (\ref{Klower}), it follows that in order for the inequality (\ref{Ksmall}) to be satisfied with $\theta_*$ given by (\ref{thetamax}), the matrix $N$ in (\ref{hlin}) has to be small enough in the sense that
\begin{equation}
\label{Nsmall}
    \|N\sqrt{P_*}\|^2
    =
    \frac{1}{2\theta_*}
    <
    \frac{1}{4K}
    \<
    -\frac{1}{2n} \Tr A .
\end{equation}
Also, the fulfillment of the sufficient condition (\ref{Nov}) for (\ref{RT}) at any $T>0$, with the system being in the nominal invariant measure according to (\ref{Xlin*}), is secured by another constraint on the matrix  $N$:
\begin{equation}
\label{NFsmall}
    \|N F\|_\infty < 1,
\end{equation}
where the transfer function $F$ is given by  (\ref{F}), and $\|\cdot\|_\infty$ is the $\cH_\infty$-norm.
The condition (\ref{NFsmall}) originates from the limit relation \cite{BV_1985,MG_1990}
\begin{align}
\nonumber
    \lim_{T\to +\infty}&
    \Big(
        \frac{1}{T}
        \ln
        \bE_*
        \re^{\frac{1}{2}\int_0^T |NX_t|^2\rd t}
    \Big)\\
\label{QEFrate}
    & =
    -
    \frac{1}{4\pi}
    \int_{\mR}
    \ln\det (I_m - N\Sigma(\omega)N^\rT)\rd \omega,
\end{align}
where $\Sigma$ is the spectral density (\ref{Sigma}), so that $N\Sigma N^\rT$ is the spectral density of the stationary Gaussian process $N X$ (computed for the case when $X$ is a stationary Gaussian process governed by the nominal SDE (\ref{Xlin*})). However, (\ref{NFsmall}) can only affect the link between $\bE (|h|^2)$ and the noise relative entropy rate on the right-hand side of (\ref{Epsi2}), whereas the inequality (\ref{dev}) and all its corollaries remain valid regardless of this link.

\section{Perturbed Langevin Dynamics Example}
\label{sec:ex}

\begin{figure}
\begin{center}
\includegraphics[width=8.4cm]{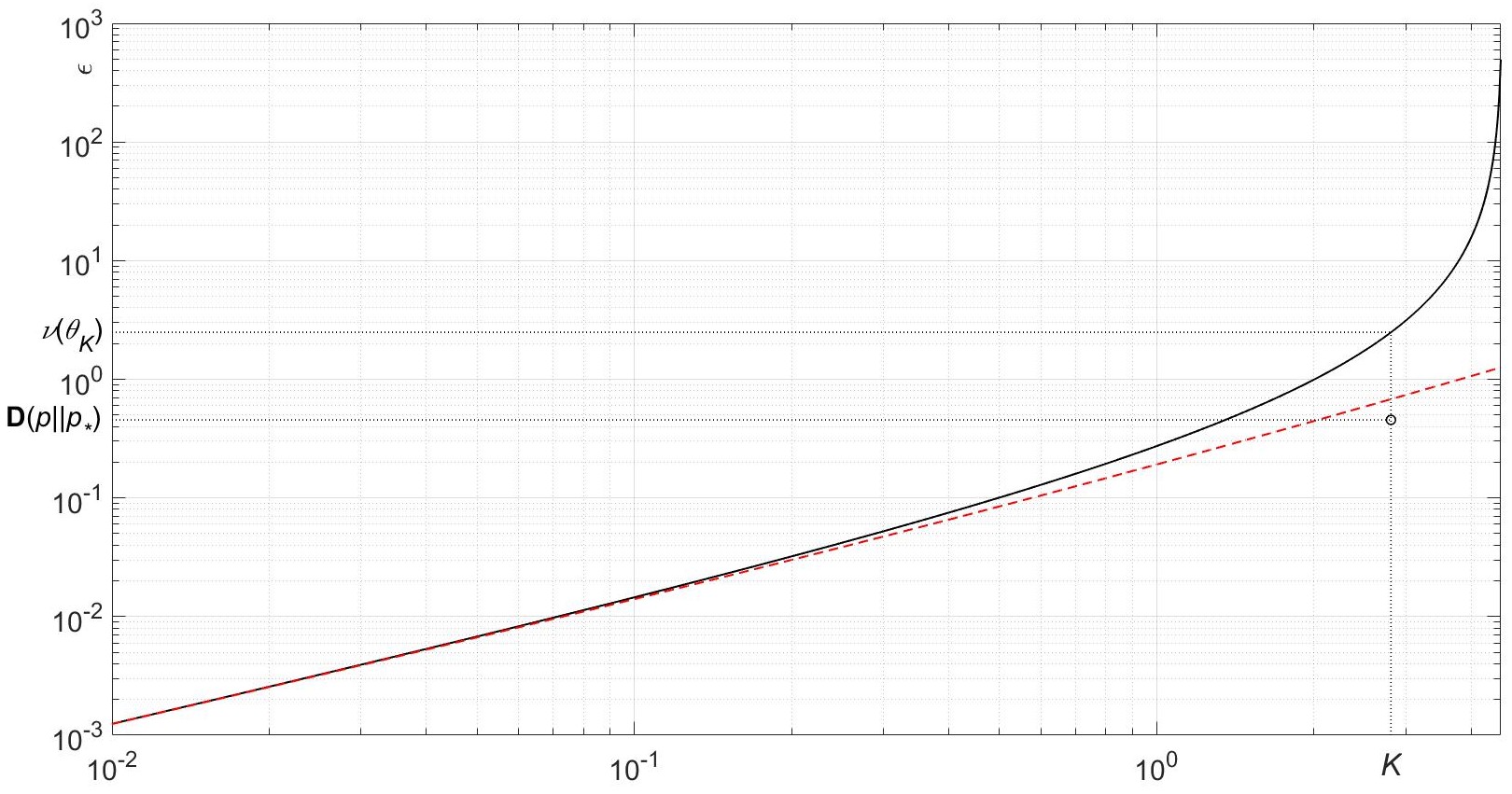}    
\caption{The relative entropy bound $\nu(\theta_K)$ from (\ref{bDupper1}) as a function (the solid black curve) of the gain coefficient $K$ from (\ref{K}),  along with its value from (\ref{lammuK}) and the exact value $\bD(p\|p_*)$ from  (\ref{bDexact}) marked by ``$\circ$''. Also graphed (the dashed red curve) is the truncated asymptotic expansion (\ref{nuasyK}).   }  \label{fig:Keps}                                 
\end{center}                                 
\end{figure}

For a comparison of the upper bound (\ref{bDupper1}) with the exact value
\begin{equation}
\label{bDexact}
    \bD(p\|p_*)
     =
    \frac{1}{2}(\Tr \chi - \ln\det \chi - n),
    \quad
    \chi:= P_*^{-1}P
\end{equation}
of the 
Kullback-Leibler relative entropy for the invariant Gaussian PDFs (\ref{p*gauss}), (\ref{pgauss}) (see, for example, \cite[Lemma~9.1]{VP_2015} and references therein), consider
a linear stochastic system whose nominal dynamics (\ref{Xlin*}) are specified by
\begin{equation}
\label{ABsigtau}
    A := -R,
    \quad
    B := \sigma I_n,
    \quad
    \sigma :=\sqrt{2\tau},
    \quad
    \tau >0,
\end{equation}
where $n=m$, and $R$ is a real positive definite symmetric matrix of order $n$. The resulting form of the SDE (\ref{Xlin*}) is
\begin{equation}
\label{Xlin*exnom}
    \rd X_t = - R X_t \rd t+\sigma \rd V_t
\end{equation}
and describes the Langevin dynamics \cite{Z_2001} for the velocity $X_t$ of a particle of unit mass in $\mR^n$ with a damping matrix $R$, or, alternatively,  a noisy gradient descent for minimising $\frac{1}{2}\|x\|_R^2$ over $x \in \mR^n$. In accordance with (\ref{DBB}),  (\ref{ABsigtau}) yields
\begin{equation}
\label{Dsig}
  D = \sigma^2 I_n = 2\tau I_n,
\end{equation}
and, by (\ref{P*ALE}), the nominal invariant covariance matrix in (\ref{P*}) is given by
\begin{equation}
\label{P*ex}
    P_* = \tau R^{-1}.
\end{equation}
The corresponding Gaussian PDF in (\ref{p*gauss}) is a Boltzmann equilibrium PDF
$    p_*(x)
    =
    \frac{\sqrt{\det R}}{(2\pi\tau)^{n/2}}
    \re^{-
        \frac{1}{2\tau}
        \|x\|_R^2}$,
$
    x \in \mR^n
$,
where $\tau$ plays the role of  an absolute temperature  parameter. Therefore, the ellipticity and logarithmic concavity constants $\lambda$, $\mu$ in (\ref{lamgauss}), (\ref{mugauss}) and the coefficient $K$ in  (\ref{K})  are
\begin{equation}
\label{lammuK}
    \lambda = \sigma^2 = 2\tau,
    \ \ \
    \mu = \frac{1}{\tau}\lambda_{\min}(R),
    \ \ \
    K = \frac{1}{2\lambda_{\min}(R)}.
\end{equation}
In the presence of a linear noise drift (\ref{hlin}), specified by a matrix $N \in \mR^{n\x n}$,   the perturbed SDE (\ref{Xlin}) takes the form
\begin{equation}
\label{Xlin*ex}
    \rd X_t = (\sigma N - R)X_t \rd t + \sigma \rd V_t,
\end{equation}
replacing (\ref{Xlin*exnom}).
Similarly to Gershgorin's circle argument \cite{HJ_2007}, if the matrix $N$ satisfies
\begin{equation}
\label{N+Nsmall}
    \lambda_{\max}(N+N^\rT)
    <
    \frac{2}{\sigma}
    \lambda_{\min}(R)
\end{equation}
(in particular, whenever $N$ is small enough in the sense that  $\|N\| < \frac{1}{\sigma}\lambda_{\min}(R)$), then
the matrix $\sigma N - R$ in (\ref{Xlin*ex}) is Hurwitz. Also note that, in view of (\ref{ABsigtau}), (\ref{Dsig}),  (\ref{P*ex}), the matrix (\ref{H}) vanishes: $H = -R +\tau (\tau R^{-1})^{-1} =0$, whereby the condition (\ref{PiH}) is satisfied for any $P$, and an equality in (\ref{PPnorm}) is indeed achieved at  a symmetric matrix $N = \frac{1}{2} \sigma\Pi$.
For a numerical example with dimensions $n = m = 4$, we used (\ref{ABsigtau}) with $\tau = 0.3463$ and the following matrix:
$$
    R
    =
    {\scriptsize\begin{bmatrix}
     1.7833  &   0.5123 &    -0.1783  &   0.1760\\
    0.5123   & 5.2275   &-3.4186      & -1.7825\\
   -0.1783   &-3.4186   &4.3321       & 0.2209\\
    0.1760   &-1.7825   &0.2209       & 1.4656
    \end{bmatrix}}
$$
(its smallest eigenvalue is $\lambda_{\min}(R)=0.1779$), so that $K = 2.8099$ in (\ref{lammuK}). The perturbed system (\ref{Xlin*ex}) was considered with
$$
    N
    =
    {\scriptsize\begin{bmatrix}
   -0.0291 &   0.0520 &  -0.0007 &  -0.0424\\
   -0.0807 &   0.0783 &   0.0474 &  -0.0066\\
    0.0570 &  -0.0590 &   0.0574 &   0.0137\\
    0.0091 &   0.0333 &   0.0425 &   0.1638
    \end{bmatrix}},
$$
which satisfies the strict inequality in (\ref{Nsmall}) with $\theta_* = 9.1946$,  (\ref{NFsmall}) with $\|NF\|_\infty = 0.7807$ (thus making (\ref{QEFrate}) valid in this example), and (\ref{N+Nsmall}).   The invariant covariance matrix $P$, found by solving (\ref{PALE}) with $A$, $B$, $D$ from (\ref{ABsigtau}), (\ref{Dsig}), is given by
$$
    P
    =
    {\scriptsize\begin{bmatrix}
    0.3949 &  -0.5799 &  -0.3971 &  -0.7857\\
   -0.5799 &   1.6852 &   1.1883 &   2.1990\\
   -0.3971 &   1.1883 &   0.9194 &   1.5359\\
   -0.7857 &   2.1990 &   1.5359 &   3.1403
   \end{bmatrix}},
$$
which, together with (\ref{P*ex}), yields the following exact value (\ref{bDexact}) for the relative entropy:
$\bD(p\|p_*)= 0.4544$. Its upper bound $\nu(\theta_K) = 2.4894$ was computed according to (\ref{bDupper1}) of Theorem~\ref{th:bDnu} using (\ref{Psi'gauss}), (\ref{nugauss}) and the parameterisation (\ref{tpar})  of (\ref{Keps}) in the $(K,\eps)$-plane shown in Fig.~\ref{fig:Keps}.
This figure also visualises the truncated asymptotic expansion for $\nu(\theta_K)$ from (\ref{nuasyK}) of Theorem~\ref{th:asy} whose coefficients were found using the Gaussian cumulants (\ref{Psi'Psi''gauss}).

\section{Conclusion}
\label{sec:conc}

We have considered a class of stochastic systems governed by SDEs whose driving noise deviates from the standard Wiener process due to a state-dependent drift. An upper bound on a diffusion-weighted Dirichlet form for the logarithmic PDF ratio (corresponding to the actual invariant measure of the perturbed system and  its nominal counterpart in the white-noise case) has been obtained in terms of the Kullback-Leibler relative entropy rate of the input noise. This bound has been shown to be achievable (and the noise drift which saturates this inequality provided) under the preservation of the PDF ratio by the divergenceless steady-state  probability flux associated with the nominal FPKE.
Together with the logarithmic Sobolev inequality, the Dirichlet form bound  has been used for obtaining (and studying the asymptotic behaviour of)  an upper bound on the relative entropy of the perturbed invariant PDF in terms of quadratic-exponential moments of the noise drift in the case of a uniformly elliptic diffusion matrix and a strongly logarithmically concave nominal invariant PDF. These results have been illustrated for perturbations of Gaussian invariant measures in linear stochastic systems with linear noise drifts, including a perturbed Langevin dynamics example.

\begin{ack}
Support from the Australian Research Council grants DP210101938, DP200102945 is gratefully acknowledged.
\end{ack}



\begin{thebibliography}{99}
\bibitem{A_1989}
V.I.Arnold, \emph{Mathematical Methods of Classical Mechanics}, 2nd Ed., Springer, New York, 1989.

\bibitem{BE_1985}
D.Bakry, and M.Emery, Diffusions hypercontractives, in ``Seminaire de
Probabilitts XIX,'', Lecture Notes in Mathematics, vol. 1123, Springer-Verlag, New York/Berlin, 1985,  pp. 175--206.
\bibitem{Beghi_1994}
A.Beghi, Continuous-time Gauss-Markov processes with fixed
reciprocal dynamics, {\it J. Math. Sys.
Estim. Contr.}, vol. 4, no. 4, 1994, pp. 1--24.
\bibitem{BFP_2002}
A.Beghi, A.Ferrante, and M.Pavon, How to steer a quantum system
over a Schr\"{o}dinger bridge, {\it Quant. Inform.
Process.}, vol. 1, no. 3, 2002, pp. 183--206.
\bibitem{BV_1985}
A.Bensoussan, and J.H. van Schuppen, Optimal control of partially observable stochastic systems with an exponential-of-integral performance index, \emph{ SIAM J. Control Optim.}, vol.  23, no. 4, 1985,  pp. 599--613.
\bibitem{Blaqiere_1992}
A.Blaqui\'{e}re, Controllability of a Fokker-Planck equation, the
Schr\"{o}dinger system, and a related stochastic optimal control
(revised version),  {\it J. Dynam. Contr.}, vol. 2,
no. 3, 1992, pp. 235--253.
\bibitem{BR_1995}
V.I.Bogachev, and M.Rockner,
Regularity of invariant measures on finite and infinite dimensional spaces and applications,
\emph{J. Funct. Anal.}, vol. 133, no. 1, 1995, 168--223.
\bibitem{BKR_1996}
V.I.Bogachev, N.V.Krylov, and M.R\"{o}ckner,
Regularity of invariant measures:
the case of non-constant diffusion part,
\emph{J. Funct. Anal.}, vol. 138, 1996, 223--242.
\bibitem{BKR_2006}
V.I.Bogachev, N.V.Krylov, and M.R\"{o}ckner,
Elliptic equations for measures:
regularity and global bounds of densities,
\emph{J. Math. Pures Appl.}, vol.  85, 2006, pp. 743--757.
\bibitem{BKRS_2015}
V.I.Bogachev, N.V.Krylov, M.R\"{o}ckner, and S.V.Sha\-posh\-ni\-kov,
\textit{Fok\-ker–Planck–Kolmogorov Equations}, 
American Mathematical Society,
Providence, Rhode Island, 2015.

\bibitem{BSV_2016}
V.I.Bogachev, S.V.Shaposhnikov, and A.Yu.Veretennikov, Differentiability of solutions of
stationary Fokker-Planck-Kolmogorov equations with respect to a parameter, \emph{Discrete
and Continuous Dynamical Systems Series A}, vol. 36, no. 7, 2016, pp. 3519--3543.

\bibitem{BKS_2023}
V.I.Bogachev, E.D.Kosov, and A.V.Shaposhnikov,
Regulari- \\ ty of solutions to Kolmogorov equations
with perturbed drifts,  \emph{Potential Anal.},  vol. 58, 2023, pp. 681--702 (published 21 September 2021).


\bibitem{B_1996}
A.Boukas, Stochastic control of operator-valued processes in boson Fock space, \emph{Russian
J. Math. Phys.}, vol. 4, 1996, pp. 139--150.

\bibitem{B_1967}
L.M.Bregman, The relaxation method of finding the common points of convex sets and its application to the solution of problems in convex programming,  \emph{USSR Comput. Math. Math. Phys.}, vol. 7, no. 3, 1967, pp. 200--217.
\bibitem{CL_1994}
P.Cattiaux,  and C.L\'{e}onard, Minimization of the Kullback information of diffusion processes,
\emph{Ann. Inst. H. Poincare}, vol.  30, no. 1, 1994, 83--132.
\bibitem{CR_2007}
 C.D.Charalambous and  F.Rezaei, Stochastic uncertain systems
subject to relative entropy constraints: induced norms and
monotonicity properties of minimax games, {\it IEEE Trans. Autom.
Contr.}, vol. 52, no. 4, 2007,  pp. 647--663.
\bibitem{CGP_2016}
Y.Chen, T.T.Georgiou, and M.Pavon, On the relation between optimal transport and Schr\"{o}dinger bridges: a stochastic control viewpoint, \emph{J. Optim. Theory Appl.}, vol.  169, 2016, pp. 671--691.
\bibitem{CT_2006}
T.M.Cover, and J.A.Thomas, {\it Elements of Information Theory},   Wiley, Hoboken, New Jersey, 2006.
\bibitem{DaiPra_1991}
P.Dai~Pra, A stochastic control approach to reciprocal diffusion
processes, {\it Appl. Math. Optim.}, vol. 23,
no. 1, 1991, pp. 313--329.
\bibitem{DE_1997}
P.Dupuis, and R.S.Ellis, {\it A Weak Convergence Approach to the
Theory of Large Deviations}, Wiley, 1997.
\bibitem{DJP_2000}
 P.Dupuis, M.R.James, and I.R.Petersen, Robust properties of
risk-sensitive control, {\it Math. Contr. Sign. Sys.}, vol.
13, 2000, pp. 318--332.
\bibitem{E_2008}
L.C.Evans, {\it Partial Differential Equations}, American Mathematical Society, Providence, Rhode Island, 2008.
\bibitem{F_1982}
W.H.Fleming, Logarithmic transformations and stochastic control,
Lecture Notes in Control and Information Sciences, vol. 42, 1982, pp.
131--141.
\bibitem{Girsanov_1960}
 I.V.Girsanov, On transforming a certain class of stochastic
processes by absolutely continuous substitution of measures, {\it
Theor. Probab. Appl.}, vol. 5, no. 3, 1960, pp. 285--301.
\bibitem{G_1975}
L.Gross,
Logarithmic Sobolev inequalities,
\emph{American Journal of Mathematics}, vol. 97, no. 4, 1975,
pp. 1061--1083.

\bibitem{G_1993}
L.Gross, Logarithmic Sobolev inequalities and contractivity properties of semigroups,  In: G.Dell'Antonio, and U.Mosco (Eds),  Dirichlet Forms,  Lecture Notes in Mathematics, vol. 1563, Springer, Berlin, Heidelberg, 1993, pp. 54--88.
\bibitem{H_1967}
L.H\"{o}rmander, Hypoelliptic second order differential equations, \textit{ Acta Math.}, vol.  119, 1967, pp.  147--171.

\bibitem{HJ_2007}
 R.A.Horn, and C.R.Johnson, {\it Matrix Analysis}, Cambridge
University Press, New York, 2007.
\bibitem{J_1973}
D.H.Jacobson, Optimal stochastic linear systems with exponential performance
criteria and their relation to deterministic differential games, \emph{IEEE Trans. Autom.
Control}, vol. 18, 1973, pp. 124--31.
\bibitem{KS_1991}
I.Karatzas, and S.E.Shreve,
\emph{Brownian Motion and Stochastic Calculus}, 2nd Ed.,
Springer, New York, 1991.
\bibitem{KS_1972}
H.Kwakernaak, and R.Sivan,
\textit{Linear Optimal Control Systems},
Wiley, New York, 1972.
\bibitem{L_1992}
M.Ledoux,
On an integral criterion for hypercontractivity
of diffusion semigroups and extremal functions,
\emph{J. Funct. Anal.}, vol. 105, 1992, 444--465.
\bibitem{LS_2001}
R.S.Liptser, and A.N.Shiryaev,
\emph{Statistics of Random Processes II:
Applications},
Springer,  Berlin, 2001.
\bibitem{MPR_2005}
G.Metafune, D.Pallara, and A.Rhandi, Global properties of invariant measures, \emph{J. Funct. Anal.}, vol.  223, 2005, pp. 396--424. 
\bibitem{Mikami_1990}
T.Mikami, Variational processes from the weak forward equation,
{\it Commun. Math. Phys.}, vol. 135, 1990, pp. 19--40.
\bibitem{MG_1990}
D.Mustafa, and K.Glover, \emph{Minimum Entropy $H_\infty$ Control}, Springer-Verlag, Berlin, 1990.
\bibitem{Nelson_2001}
E.Nelson, {\it Dynamical Theories of Brownian Motion}, 2nd Ed.,
Princeton University Press, 2001.
\bibitem{N_1973}
A.A.Novikov, On an identity for stochastic integrals,
\emph{Theor. Probab. Appl.},  vol. 17, no. 4, 1973, pp. 717--720.
\bibitem{O_2000}
B.{\O}ksendal, \emph{Stochastic Differential Equations}, 5th Ed., Springer-Verlag, Berlin, 2000.
\bibitem{PF_2013}
M.Pavon,  and A.Ferrante, On the geometry of maximum entropy problems, \emph{SIAM Rev.}, vol. 55, no. 3, 2013, pp. 415--439.
\bibitem{PUS_2000}
I.R.Petersen, V.A.Ugrinovskii, and A.V.Savkin, {\it Robust Control
Design Using $H^{\infty}$ Methods}, Springer, London, 2000.
\bibitem{P_2006}
I.R.Petersen, Minimax LQG control, \emph{Int. J. Appl. Math. Comput. Sci.}, vol. 16, 2006, pp. 309--323.
\bibitem{R_1996}
H.Risken,
\emph{The Fokker-Planck Equation: Methods of Solution and Applications}, 2nd Ed., Springer, Berlin, 1996. 
\bibitem{S_1959}
A.J.Stam,
Some inequalities satisfied by the quantities
of information of Fisher and Shannon,
\emph{Information and Control}, vol.  2, no. 2, 1959, 101--112.

\bibitem{S_2008}
D.W.Stroock,
\textit{Partial Differential Equations for Probabilists},
Cambridge University Press, Cambridge, 2008.
\bibitem{TEM_2018}
T.Tanaka, P.M.Esfahani,  and S.K.Mitter, LQG control with minimum directed information: semidefinite programming approach, \emph{IEEE Trans. Automat. Contr.}, vol. 63, no. 1, 2018, pp. 37--52.
\bibitem{T_2009}
A.B.Tsybakov, \emph{Introduction to Nonparametric Estimation}, Springer, New York,
2009.
\bibitem{UP_2001}
 V.A.Ugrinovskii and I.R.Petersen, Minimax LQG
control of stochastic partially observed uncertain systems, {\it
SIAM J. Contr. Optim.}, vol. 40, no. 4, 2001, pp. 1189--1226.
\bibitem{VKS_1995}
I.G.Vladimirov, A.P.Kurdjukov, and A.V.Semyonov, Anisotropy of signals and entropy of linear time-invariant systems, \emph{Doklady Mathematics}, vol. 342, no. 5, 1995, pp. 583--585 (Russian) (English translation in vol. 51, no. 3, pp. 388--390).
\bibitem{VP_2010}
I.G.Vladimirov, and I.R.Petersen,
Minimum relative entropy state transitions in linear stochastic
systems: the continuous time case, In:
Proc. 19th International Symposium on Mathematical Theory of Networks and Systems (MTNS), 5-9 July 2010, Budapest, Hungary, pp. 51--58.
\bibitem{VP_2015}
I.G.Vladimirov, and I.R.Petersen,     State distributions and
minimum relative entropy noise sequences
    in uncertain stochastic systems:
    the discrete time case,
\emph{SIAM J. Contr. Optim.}, vol. 53, no. 3, 2015, pp. 1107--1153.
\bibitem{VPJ_2018}
I.G.Vladimirov, I.R.Petersen, and M.R.James, Multi-point Gaussian states, quadratic-exponential cost functionals, and large deviations estimates for linear quantum stochastic systems, \textit{Appl. Math. Optim.}, vol. 83, no. 1, 2021, pp. 83--137 (published online 24 July 2018).

\bibitem{VPJ_2021}
I.G.Vladimirov, I.R.Petersen, and M.R.James, Quadratic-exponential functionals of Gaussian quantum processes, \emph{Infinite-dimensional Analysis, Quantum Probability and Related Topics}, vol. 24, no. 4, 2021, pp. 2150024.


\bibitem{V_2022}
I.G.Vladimirov,
Probabilistic bounds with quadratic-exponential moments for
quantum stochastic systems, IFAC World Congress,  9-14 July 2023, Yokohama, Japan, accepted (preprint: 	arXiv:2211.12161 [quant-ph], 22 November 2022).

\bibitem{W_1981}
P.Whittle, Risk-sensitive linear/quadratic/Gaussian control, \emph{Adv. Appl. Probab.},
vol. 13, 1981, pp. 764--77.
\bibitem{YKT_2023}
A.V.Yurchenkov, A.Yu.Kustov, and V.N.Timin,
The sensor network estimation with dropouts: anisotropy-based approach,
\emph{Automatica}, 
vol. 151, 2023, pp. 110924.

\bibitem{Z_2001}
R.Zwanzig, \textit{Nonequilibrium Statistical Mechanics}, Oxford University Press, New
York, 2001.


\end{thebibliography}
\end{document}